\documentclass[letterpaper]{article}

\usepackage{arxiv}

\usepackage{ART-SFVEM}

\usepackage{lineno,hyperref}

\usepackage{amsmath, amssymb, latexsym, multicol, amsthm, bm, placeins}
\usepackage[export]{adjustbox}
\usepackage{color}
\usepackage{multirow}

\makeatletter\@addtoreset{equation}{section}\makeatother

\renewcommand{\shorttitle}{A residual a posteriori error estimate for the Stabilization-free Virtual Element
  Method}

\numberwithin{equation}{section}

\makeatletter
\renewcommand{\@biblabel}[1]{#1\hfill \hspace{0.2cm}}
\makeatother

\begin{document}
\title{A residual a posteriori error estimate for the Stabilization-free Virtual Element
  Method}  
\author{
  Stefano Berrone \thanks{Dipartimento di Scienze Matematiche “Giuseppe Luigi Lagrange”, Politecnico di Torino, Corso Duca degli Abruzzi 24, Torino, 10129, Italy. stefano.berrone@polito.it (S. Berrone), andrea.borio@polito.it (A. Borio), davide.fassino@polito.it (D. Fassino), francesca.marcon@polito.it (F. Marcon).}
  \And
  Andrea Borio \footnotemark[1]
  \And
  Davide Fassino \footnotemark[1]
  \And
  Francesca Marcon \footnotemark[1]
  }
%\newcommand{\poliaddress}{Dipartimento di Scienze
%  Matematiche, Politecnico di Torino,
%  % Corso Duca degli Abruzzi 24   
%  Torino, 10129, Italy
%} 

%\author{
%  Stefano Berrone \thanks{\poliaddress (stefano.berrone@polito.it)}
%  \And
%  Andrea Borio \thanks{\poliaddress (andrea.borio@polito.it)}
%  \And
%  Davide Fassino \thanks{\poliaddress (davide.fassino@polito.it)}
%  \And
%  Francesca Marcon \thanks{\poliaddress (francesca.marcon@polito.it)}
%  }
\maketitle

\begin{abstract}
In this work, we present the a posteriori error analysis of Stabilization-Free
Virtual Element Methods for the 2D Poisson equation. The abscence of a
stabilizing bilinear form in the scheme allows to prove the equivalence between
a suitably defined error measure and standard residual error estimators, which
is not obtained in general for stabilized virtual elements. Several numerical
experiments are carried out, confirming the expected behaviour of the estimator
in the presence of different mesh types, and robustness with respect to jumps of
the diffusion term.

\end{abstract}
\keywords{Virtual Element Method, A posteriori error analysis, Stabilization-free.
	}

\section{Introduction}
Virtual Element Methods (VEM) are a powerful technology enabling the solution to
partial differential equations on general polytopal meshes. First introduced in
\cite{Beirao2013a,Beirao2013b}, this family of Galerkin methods is based on the
definition of discrete spaces whose basis functions are not everywhere known
analytically. Suitable discrete bilinear forms can be defined from the
degrees of freedom, exploiting the computability of suitably defined polynomial projections
of basis functions, which in general have a non-trivial kernel. Classical VEM
bilinear forms also include a stabilizing operator, defined to take care of the
kernel of the polynomial projections involved. Recently, the scientific
community has invested a considerable effort in studying polynomial
projection operators that are computable from the degrees of freedom of a VEM
function and also stable. Starting from the seminal works
\cite{BBM_ComparisonStabFree,DAltri2021,Berrone2023AfirstOrderStabFree}, novel
stabilization-free VEM have been devised for elliptic problems, based on
polynomial projections of the gradients of basis functions on high degree
polynomial spaces. In
\cite{berrone_lowest_2025,Berrone2024a,BBFM2025,Borio2024,Borio2024b,Marcon2025}, the
theoretical foundations of the method were developed, considering the primal and
mixed formulations of elliptic problems. Moreover, these discretization
strategies have been applied to various problems of interest in computational
mechanics \cite{Xu2024a, Xu2023, Chen2023a, Chen2023b, Xu2024b, Lamperti2022,
  LIGUORI2024117281}.

This paper deals with the \emph{a posteriori} error estimation of the scheme
presented in \cite{Berrone2024a,BBFM2025}. In the context of classical
stabilized VEM, the development of a posteriori error estimates can be found in
\cite{berrone_residual_2017, Cangiani2017, ABBDVW, MUNAR2024182} and adaptive
schemes exploiting these results have been succesfully devised, addressing the
non-trivial issue of defining quality-preserving refinement strategies for
polygons \cite{berrone_reliable_2019, berrone_refinement_2021,
  berrone_effective_2025}. However, in the a posteriori analysis of standard VEM
the stabilization operator is found to be an issue in deriving the required
equivalence between the error and the error estimator. Recently, some
stabilization-free a posteriori error bounds were proved for stabilized VEM
\cite{Beirao2023b,Canuto2023,Berrone2024b}, though these results are limited to
certain classes of meshes. In this work, the abscence of a stabilizing bilinear
form in stabilization-free schemes enables to prove rigorously the equivalence
between a suitably defined error measure and classical residual error
estimators. A campaign of numerical tests confirms that the ratio between the
error and the error estimator is indeed asymptotically constant. To streamline
the analysis, estimates are presented in the case of a piecewise constant
diffusivity coefficient $\K$. 
The extension to variable $\K$ would exploit analogous residual estimators and 
would  yield a dependence of the estimate on local oscillations of the diffusivity 
coefficient, similarly to what happens in the \emph{a posteriori} analysis of finite element methods.

The outline of the paper is as follows: in Sections \ref{sec:ModelProblem} and
\ref{sec:DiscreteProblem} we introduce the model problem and describe the
stabilization-free Virtual Element Method used to discretize the problem. In
Section \ref{sec:apost} we present the \emph{a posteriori} error analysis of the
numerical method. Section \ref{sec:NumericalResults} is devoted to numerical
results on a squared and an L-shaped domain. Finally, in Section
\ref{sec:conclusions} we summarize the results and draw conclusions.
\section{Model problem}\label{sec:ModelProblem}
Let $\Omega \subset \mathbb{R}^2$ be a bounded open set. We are interested in solving the following Poisson problem:
\begin{equation}
  \label{eq:modelProblem}
  \begin{cases}
    -\nabla \cdot\left( \K \,\nabla u \right) = f  & \text{ in $\Omega$,}\\
       u= 0 &\text{ on $\partial \Omega$,}
  \end{cases}
\end{equation}
where the loading term $f\in L^2(\Omega)$ and $\K$ is a symmetric tensor $\K\in(L^\infty(\Omega))^{2\times 2}$ satisfing:
\begin{equation}\label{eq:constraintsK}
   \K_0|\bm{v}|^2  \le \bm{v}\cdot \K(\bm{x})\bm{v}\le \K_1 |\bm{v}|^2, \quad \forall \bm{v}\in \mathbb{R}^2, \quad \text{for a.e. }\,\bm{x}\in\Omega,
\end{equation}
where $\K_0$ and $\K_1$  are positive constants and $\abs{\cdot}$ denotes the euclidean norm.
We define the bilinear form $\a{\cdot}{\cdot}\colon \sobh[0]{1}{\Omega}\times \sobh[0]{1}{\Omega} \to \mathbb{R}$ as 
\begin{equation}
  \label{eq:defBilFormA}
  \a{u}{v} := \int_\Omega\left(\K \,\nabla u\right) \cdot\nabla v,\quad \forall \,u,v \in \sobh[0]{1}{\Omega}.
\end{equation}
The variational formulation of \eqref{eq:modelProblem} reads as:
find $u\in \sobh[0]{1}{\Omega}$ such that,
\begin{equation}
  \label{eq:contVarForm}
  \a{u}{v} = \scal[\Omega]{f}{v} \quad \forall \,v\in\sobh[0]{1}{\Omega}, \,
\end{equation}
where $\scal[\mathcal{O}]{\cdot}{\cdot}$ denotes the scalar product in $L^2(\mathcal{O})$, where $\mathcal{O}\subseteq \Omega$.
\section{Mesh and discretization}\label{sec:DiscreteProblem}
We introduce a conforming polygonal tessellation $\Mh$ of $\Omega$. Let $E$
denote a generic polygon of $\Mh$, $h_E$ the diameter of $E$ and the mesh size
$h\coloneqq \max_{E\in\Mh}h_E$. Furthermore, let $\Xh$ denote the set of vertices of $\Mh$. 
Moreover, let $\Eh$ be the set of edges of $\Mh$, $\Ehint$ the set of edges of $\Mh$ internal to $\Omega$ and
for each $E\in\Mh$ $\Eh[E]\coloneqq\{e\in\Eh : e\subset \partial E\}$.
We assume that $\Mh$ satisfies the standard
mesh assumptions for VEM, described below.
\begin{assumption}[Mesh assumptions] \label{hp:mesh_assumption_VEM}
$\exists \,\kappa > 0$ such that
$\forall E\in\Mh$, 
\begin{itemize}
\item $E$ is star-shaped with respect to a ball of radius $\rho_E\geq \kappa h_E$,% where $h_E$ is the diameter of $E$,
\item for every edge $e\subset\partial E$, $\abs{e}=: h_e \geq \kappa h_E$.
\end{itemize}
Notice that the above conditions imply that, denoting by $\NVE$ the number of vertices of $E$, the number of vertices of each polygon $E$ has an upper bound, i.e. 
\begin{equation}\label{eq:numVerticesBounded}
\exists \NVmax>0\colon \forall E\in\Mh,\, \NVE \leq \NVmax\,.
\end{equation}
Furthermore, these conditions ensure  that the number of elements in the mesh sharing a vertex $\bs{x}\in\Xh$, denoted by $\NVEE$, remains bounded independently of $h$, i.e.
\begin{equation}\label{eq:numElementperVertexBounded}
\exists \NVEmax>0\colon \forall \bm{x}\in\Xh,\, \NVEE \leq \NVEmax\,.
\end{equation}
\end{assumption}
From now on, for simplicity of exposition, we assume $\K$ to be a constant on
each polygon $E\in\Mh$. Let $\K_E \in \mathbb{R}$ denote its value on each $E\in\Mh$.

\subsection{Discretization}

Let $\Poly{k}{\mathcal{O}}$ be the set of polynomials
of degree up to $k$ defined on an open connected set $\mathcal{O}$, and let
$\Poly{m}{\mathcal{O}} \setminus \Poly{n-1}{\mathcal{O}}$ be the set of polynomials from
degree $n$ up to $m$. For any given $E\in\Mh$, let
$\proj[\nabla]{k}{E}:\sobh{1}{E}\rightarrow\Poly{k}{E}$ be the $\sobh{1}{}$
orthogonal projection operator such that, $\forall \,u\in\sobh{1}{E} $,
\begin{equation}\label{eq:higher_order_PiNablaorthogonalitycondition}
  \begin{cases}
    \scal[E]{\nabla\left(\proj[\nabla]{k}{E} u -u \right)}{\nabla p}=0
    & \forall p\in\Poly{k}{E}
    \\
    \int_{E} u - \proj[\nabla]{k}{E} u = 0 & \text{if $k>1$}
    \\
    \int_{\partial E} u - \proj[\nabla]{k}{E} u = 0 & \text{if $k=1$}
  \end{cases}
\end{equation}
Moreover, we define the local Virtual Space of order $k$
 as
\begin{multline}\label{def:Vemspace}
  \V[E]{k}:=\{ v_h\in\sobh{1}{E} : \Delta v_h\in \Poly{k}{E}, \;\; \trace{v_h}{e}{}\in \Poly{k}{e}\;\forall e \subset \partial E, \\ v_h\in\cont{\partial E},\;\; \scal[E]{v_h-\proj[\nabla]{k}{E} v_h}{p}=0\; \forall p\in\Poly{k}{E}\setminus \Poly{k-2}{E} \}
\end{multline}
Given $v_h\in\V[E]{k}$, the chosen degrees of freedom of this space are:
\begin{itemize}
\item $\NVE$ pointwise values of $v_h$ at the vertices of the polygon,
\item if $k>1$, $k-1$ pointwise values of $v_h$ at Gauss-Lobatto quadrature
  points internal to each edge,
\item if $k>1$, $\frac{k(k-1)}{2}$ internal moments
  $\frac{1}{E}\scal[E]{v_h}{m_i}$, $\forall \, i=1,\ldots, n_{k-2}$, where
  $n_{k-2}:=\dim\Poly{k-2}{E}$ and $\{m_i\}_{i=1}^{n_{k-2}}$ is a scaled monomial basis of $\Poly{k-2}{E}$.
\end{itemize}
Then, we define the global discrete space as
$ \V{k} := \{v\in\sobh[0]{1}{\Omega}\colon v_{|E} \in \V[E]{k}\} $.

To discretize \eqref{eq:contVarForm}, we follow the approach described in
\cite{BBFM2025}.  For any given $E\in\Mh$, let $\ell_E\geq 0$ be a given natural
number and let $\mathcal{P}_{k,\ell_E}(E)$ be the polynomial space given by
\begin{equation}
  \label{eq:firstDecompositionPkell}
  \begin{aligned}
    \mathcal{P}_{k,\ell_E}
    &:= \PolyDouble{k-1}{E}\oplus \curl\left(\Poly{k+\ell_E}{E}\setminus\Poly{k}{E}\right)
    \\
    &=\bs{x}\Poly{k-2}{E}\oplus\curl\Poly{k+\ell_E}{E}\,,
  \end{aligned}
\end{equation} 
where for any $p\in\Poly{k+\ell_E}{E}$,
$\curl p =\left(\frac{\partial p}{\partial x_2},-\frac{\partial p}{\partial
    x_1}\right)$.
Then, let $\projH \nabla : \sobh{1}{E} \to \mathcal{P}_{k,\ell_E}$ be the
$\lebl{}$-projection operator of the gradient of functions in $\sobh{1}{E}$,
defined, $\forall u\in \sobh{1}{E}$, by the orthogonality condition
\begin{equation} \label{eq:def-projH} \scal[E]{\projH \nabla u}{\boldsymbol{p}}
  = \scal[E]{\nabla u}{\boldsymbol{p}} \;\; \forall \boldsymbol{p} \in
  \mathcal{P}_{k,\ell_E} \,.
\end{equation}
\begin{remark}
  Notice that for each function $u_h\in\V[E]{k}$, the projection
  $\projH \nabla u_h$ is computable exploiting only the degrees of freedom of
  $u_h$.
\end{remark}
Next, we define the discrete bilinear form $\ahE{}{}\colon \V[E]{k}\times \V[E]{k} \to \mathbb{R}$ such that
\begin{equation}\label{eq:higher_order_defahE}
  \ahE{u_h}{v_h} :=  \scal[E]{\K_E \projH \nabla u_h}{\projH \nabla v_h}, \quad \forall u_h,v_h \in \V[E]{k}\,,
\end{equation}
Summing up over all the elements of $\Mh$, we define
$\ah{}{}\colon \V{k}\times \V{k} \to \mathbb{R}$ as
\begin{equation}
  \label{eq:higher_order_defah}
  \ah{u_h}{v_h}:=\sum_{E\in\Mh}\ahE{u_h}{v_h} \quad \forall u_h,v_h \in \V{k} \,.
\end{equation}
We can state the discrete problem as: find $u_h\in\V{k}$ such that
\begin{equation}\label{eq:higher_order_discrVarForm}
  \ah{u_h}{v_h} = \sum_{E\in\Mh}\scal[E]{f_h}{v_h} \quad \forall v_h \in \V{k}\,,
\end{equation}
where $f_h \coloneqq \proj{k}{E} f$ and $\forall E\in\Mh$,
$\proj{k}{E}\colon \lebl{E} \to \Poly{k}{E}$ is the $\lebl{}$-projection,
defined $\forall u\in\lebl{E}$ by
\begin{equation}\label{eq:PiK}
  \scal[E]{\proj{k}{E} u}{p} =\scal[E]{u}{p} \quad \forall p\in\Poly{k}{E} \,.
\end{equation}
The well-posedness of \eqref{eq:higher_order_discrVarForm} is discussed in
\cite{BBFM2025}. We summarize it in the following result.

\begin{assumption}\label{hp:lE_assumption}
  We assume $\ell_E$ to be the smallest integer such that any polynomial
  $q\in \mathbb{P}_{k+\ell_E}(E)$ can be identified by a set of degrees of
  freedom which contains $k\NVE -1$ distinct moments
  $\frac{1}{\abs{\partial E}}\scal[\partial E]{q}{\pi_i}$ with respect to a
  scaled polynomial basis of the space
  $\Poly[0]{k-1}{\partial E}:=\{\pi: \pi_{|_e} \in\Poly{k-1}{e}\,, \forall
  e\subset\partial E\,, \int_{\partial E} \pi = 0 \}$.
\end{assumption}

\begin{theorem}
  Under the mesh Assumptions \ref{hp:mesh_assumption_VEM} and Assumption
  \ref{hp:lE_assumption}, there exists for any $E\in\Mh$ a constant
  $\alpha^E_\ast$ independent of $h_E$ such that
  \begin{align*}
    \norm[0,E] {\projH \nabla v_h}^2\ge \alpha^E_\ast  \norm[0,E] { \nabla v_h}^2,\qquad \qquad\forall v_h \in\V{E}.
  \end{align*}
\end{theorem}
Using the above result, it is immediate to see that
\eqref{eq:higher_order_discrVarForm} admits a unique solution under Assumption
\ref{hp:lE_assumption}. In \cite{BBFM2025} it is also proved that Assumption
\ref{hp:lE_assumption} is also a necessary condition for well-posedness.

\section{A posteriori error analysis}
\label{sec:apost}
In this work, we consider as a measure of the error between the solution 
to problem \eqref{eq:contVarForm} $u$ and the solution to problem \eqref{eq:higher_order_discrVarForm} $u_h$,
the quantity
\begin{equation}\label{def:error}
\begin{split}
    \ennorm[\K,1,\omega]{u-u_h} & \coloneqq   \sup_{w\in\sobho{\omega}{}}\frac{\sum_{ E\in \omega} \aE{u}{w} - \ahE{u_h}{w} }{\norm[0,\omega]{\sqrt{\K}\nabla w}} \\
    &=\sup_{w\in\sobho{\omega}{}}\frac{\sum_{ E\in \omega} \scal[E]{\K_E  \left( \nabla u - \projH \nabla u_h\right) }{\nabla w}}{\norm[0,\omega]{\sqrt{\K}\nabla w}},
    \end{split}
\end{equation}
where $\aE{\cdot}{\cdot}\coloneqq\a{\cdot}{\cdot}|_E$, and $\omega\subseteq \Mh$ is the union of the closure of a set of elements in $\Mh$.
\begin{remark}
    We notice that if $u|_E\in \Poly{k}{E}$ $\forall E\in \Mh$, it results that \\$ \ennorm[\K,1,\omega]{u-u_h}~=~0$. 
\end{remark}
\subsection{Preliminary results}
In this section, we discuss some preliminary results that will be used in the definition of the upper and lower bounds of the error. 
% Firstly, we report \cite[Equation 3.15]{berrone_residual_2017} as the following Lemma.
\begin{lemma}\label{lemma:projector} Let $E\in \Mh$, there exists a constant  $C_p>0$ independent of $h_E$, such that
    \begin{align*}
       \norm[0,E]{v - \proj{k}{E}v}\le C_p \frac{h_E}{\sqrt{\K_E}}\norm[0,E]{\sqrt{\K}\nabla v}, \quad \forall v\in\sobh{1}{E}.
    \end{align*}
\end{lemma}
\begin{proof} For any $E\in \Mh$, and function $v\in \sobh{1}{E}$, we employ the Poincar\'{e} inequality, i.e.
\begin{align*}
         \norm[0,E]{v - \proj{k}{E}v}\le \norm[0,E]{v - \proj{0}{E}v}\le C_p h_E\norm[0,E]{\nabla v}= C_p \frac{h_E}{\sqrt{\K_E}}\norm[0,E]{\sqrt{\K}\nabla v}.
\end{align*}
\end{proof}
    
In \cite[Theorem 11]{Cangiani2017} it has been proved the existence of the Clément quasi-interpolation operator for the VEM space, in order to give it, we recall some useful definitions.
\begin{definition} \label{def:mesh_variables} For any element $E\in\Mh$, we define
  $\tilde{\omega}_E$ as the patch of elements with non-empty intersection with $E$, i.e. sharing at
  least one vertex.  Furthermore, for any given $e\in\Eh$, let
  $\omega_e=\displaystyle\bigcup_{E\in\Mh:e\in\Eh[E]} E$.
  Furthermore, 
  ${\tilde{\omega}_e}=\bigcup_{E\in\omega_e} \tilde{\omega}_E$, and
  ${\omega_E}=\bigcup_{e\in\Eh[E]} \omega_e$.
\end{definition}
    It can be proved that 
    \begin{equation}
       \Nomegamax =  \max_{E\in \Mh}\{\# \tilde{\omega}_E\} \le \NVmax \left(\NVEmax -2\right) +1,
    \end{equation}
    where $\NVmax$ and $\NVEmax$ are defined in  
    Assumption \ref{hp:mesh_assumption_VEM} and $\# \tilde{\omega}_E$ is the cardinality of $\tilde{\omega}_E$. 
    Notice that, if $\Mh$ is a mesh of squares and $E$ is an internal element, we get
    $\# \tilde{\omega}_E = 9$, achieving the equality in the above bound.  
\begin{lemma}[Cl\'ement interpolation estimates] \label{lemma:higher_order_interpolationError}
Under Assumption \ref{hp:mesh_assumption_VEM}, for any $v\in\sobh[0]{1}{\Omega}$,  there exists a Cl\'ement quasi-interpolation operator $v_I\in\V{k}$ satisfying for each $E\in\Mh$
\begin{equation}\label{eq:BoundClement}
    \norm[0, E]{v-v_I} + h_E \norm[0, E]{\nabla(v-v_I)}\le C_{I} h_E \norm[0, \tilde{\omega}_E]{\nabla v}\,,
\end{equation}
where $C_I$ is a positive constant, depending only on the polynomial degree $k$ and on the mesh regularity. 
\end{lemma}
\begin{definition}
    For any given $e\in\Eh$, we define $\K_{\omega_e} = \sum_{E\in\omega_e}\K_E$.
\end{definition}
\begin{lemma}\label{lemma:ClementQuasiInterpolationOperator} Let $v_I$ be the Clément quasi-interpolation operator as in Lemma \ref{lemma:higher_order_interpolationError}, then there exist three positive constants $C_{I,1}$, $C_{\K,E}$ and $C_{\K_{\omega_e}}$ such that, $\forall v \in \sobho{\Omega}{}$,
\begin{align}
        \norm[0,E]{v-v_I}&\le C_{I} C_{\K, E}\frac{h_E} {\sqrt{\K_E}}\norm[0,\tilde{\omega}_E]{  \sqrt{\K} \nabla v} 
        &\forall E\in\Mh,\label{eq:BoundClementElementK} \\
     \norm[0,e]{v- v_I}&\le C_{I,1} C_{\K_{\omega_e}} \frac{h_e^{1/2}}{\sqrt{\K_{\omega_e}}}\norm[0, \tilde{\omega}_e]{\sqrt{\K}\nabla v} &\forall e \in \Eh,\label{eq:BoundClementSqrtEdgeK}
\end{align}
where $ C_{\K, E}, C_{\K_{\omega_e}}, C_{I,1}$ are positive constants, depending only on the polynomial degree $k$, on the mesh regularity and on the diffusivity coefficients $\K_E$. 
In particular,  $C_{\K, E}~\coloneqq~\sqrt{\frac{\K_E}{\min_{E'\in \tilde{\omega}_E}\{\K_{E'}\}}}$, and $C_{\K_{\omega_e}}~\coloneqq~  \sum_{E\in\omega_e}\{C^2_{\K, E}\}$.
\end{lemma}
\begin{proof} 
Considering, \eqref{eq:BoundClementElementK}, let $E\in\Mh$. From \eqref{eq:BoundClement}, we have that
\begin{align*}
    &\norm[0,E]{v-v_I} \le C_{I} h_E\norm[0, \tilde{\omega}_E]{\nabla v} \le C_{I}\frac{h_E}{\sqrt{\min_{E'\in \tilde{\omega}_E}\{\K_{E'}\}}}\norm[0, \tilde{\omega}_E]{\sqrt{\K}\nabla v}
    \\
    &\qquad= C_{I}\frac{\sqrt{\K_E}}{\sqrt{\min_{E'\in \tilde{\omega}_E}\{\K_{E'}\}}}\frac{h_E}{\sqrt{\K_E}}\norm[0, \tilde{\omega}_E]{\sqrt{\K}\nabla v} = C_{I} C_{\K, E}\frac{h_E}{\sqrt{\K_E}}\norm[0, \tilde{\omega}_E]{\sqrt{\K}\nabla v} \,.
\end{align*}
    %where $C_{\K, E}\coloneqq\sqrt{ \frac{\Kwomegae}{\Kvomegae}}$.
Then, let $e\in\Eh$ be fixed. For any $E\in\omega_e$, we recall the following scaled trace inequality:
\begin{equation}
\label{eq:scaled-trace-inequality}
  \norm[e]{v}^2\le C_{tr}\left(h^{-1}_e\norm[E]{v}^2 + h_e\norm[E]{\nabla v}^2\right),  \quad \forall v\in H^1(E). 
\end{equation}
Then, by  \eqref{eq:scaled-trace-inequality},  \eqref{eq:BoundClement}, and Assumption \ref{hp:mesh_assumption_VEM}, we have, $\forall E\in\omega_e$,
\begin{align*}
    \norm[0,e]{v-v_I}^2&\le C_{tr}\left(h^{-1}_e\norm[E]{v-v_I}^2 + h_e\norm[E]{\nabla ( v -v_I)}^2\right)
    \\
    &\le C_{tr}\left( C_{I}^2 h^{-1}_e h^2_E\norm[0, \tilde{\omega}_E]{\nabla v}^2 + C_{I}^2 h_e\norm[0, \tilde{\omega}_E]{\nabla v}^2\right)
    \\
    &= C_{tr}C_{I}^2 \left( \frac{h^2_E}{h^2_e} + 1 \right)h_e\norm[0, \tilde{\omega}_E]{\nabla v}^2
    \\
     &\le C_{tr}C_{I}^2\left( \frac{1}{\kappa^2} + 1 \right)C^2_{\K, E}\frac{h_e}{\K_E}\norm[0, \tilde{\omega}_E]{\sqrt{\K}\nabla v}^2 \,.
     % \\
     % &\le C_{I,1}^2 C^2_{\K, E}\frac{h_e}{\K_E}\norm[0, \tilde{\omega}_E]{\sqrt{\K} \nabla v}^2
\end{align*}
Hence, 
\begin{equation*}
     \K_E \norm[0,e]{v-v_I}^2 \leq  C_{tr}C_{I}^2\left( \frac{1}{\kappa^2} + 1 \right)C^2_{\K, E}h_e\norm[0, \tilde{\omega}_E]{\sqrt{\K}\nabla v}^2 
     \quad \forall E\in\omega_e\,.
\end{equation*}
 Summing up over $E\in\omega_e$, we obtain 
\begin{equation*}
    \K_{\omega_e}\norm[0,e]{v-v_I}^2 \le C_{I,1}^2 C_{\K_{\omega_e}}^2 h_e\norm[0, \tilde{\omega}_{e}]{\sqrt{\K}\nabla v}^2 \,.
\end{equation*}
where $ C_{I,1}^2\coloneqq  C_{tr}C_{I}^2\left( \frac{1}{\kappa^2} + 1 \right)$.
\end{proof}

\begin{lemma}[Galerkin orthogonality]\label{lemma:GalerkinOrthogonality}
  Let $u$ be the solution to the continuous problem \eqref{eq:contVarForm} and $u_h \in\V{k}$ the
  solution to the discrete problem \eqref{eq:higher_order_discrVarForm}, it holds that
\begin{align*}
 \a{u}{w_h}- \ah{u_h}{w_h}= \scal[\Omega]{f-f_h}{ w_h} \quad \forall w_h \in \V{k}.
\end{align*}
\end{lemma}
\begin{proof}
It follows immediately from the definition of the continuous problem\eqref{eq:contVarForm} and the discrete problem \eqref{eq:higher_order_discrVarForm}.
\end{proof}
\subsection{Upper bound of the error}
\begin{theorem}[Upper bound]\label{theorem:upperBound}
  Let $u$ be the solution to the continuous problem \eqref{eq:contVarForm} and $u_h \in\V{k}$ be the solution to the discrete problem \eqref{eq:higher_order_discrVarForm}. Then, there exists a
  constant $C_U>0$ independent of
  $h_E$, depending on the mesh regularity, 
  such that
  \begin{align}
    \ennorm[\K,1,\Mh]{u-u_h}^2  \le C_U  \, \sum_{E\in \Mh} \left(\eta^2_E +\mathcal{F}^2_E \right),
  \end{align}
  where $\ennorm[\K,1,\Mh]{ \cdot }$ is defined by \eqref{def:error} and
  \begin{align}
    &\eta^2_E \coloneqq \frac{h^2_E}{\K_E}\norm[0, E]{r_E}^2+ \frac{1}{2}\sum_{e\in \Eh[E]\cap\Ehint}\frac{h_e}{\K_{\omega_e}} \norm[0,e]{j_e}^2,\\
    &\mathcal{F}^2_E \coloneqq \frac{h^2_E}{\K_E} \norm[0,E]{f-f_h}^2 \,,
  \end{align}
  and
  \begin{align}
    \label{eq:defRE}
    r_E &\coloneqq f_h + \divt \left(\K_E \, \projH \nabla u_h\right),
    \\
    \label{eq:defJE}
    j_e &\coloneqq [[\K \, \projHOmega \nabla u_h]]_e :=
          \sum_{E\in\omega_e}
          \K_{E}  \, \projH\nabla u_h \cdot \bm{n}^e_E \,,
         % \left(\K_{E_1}  \, \projHOmega\nabla u_h\right)|_{E_1} \cdot \bm{n}_1 + \left(\K_{E_2}  \, \projHOmega \nabla u_h\right)|_{E_2} \cdot \bm{n}_2,
  \end{align}
  where $\bm{n}^e_E$ is the normal vector to $e$ pointing outward with respect to $E\in\omega_e$.
  % $\{E_i\}_{i=1,2}$ and $\K_{\omega_e} = \K_{E_1} +\K_{E_2}.$
\end{theorem}
\begin{proof}
Let $w\in\sobho{\Omega}{}$ and $w_I\in\V{k}$ such that it satisfies Lemma \ref{lemma:higher_order_interpolationError}. By using Lemma \ref{lemma:GalerkinOrthogonality} we have that
\begin{align*}
&\sum_{E\in \Mh} \scal[E]{\K_E \left( \nabla u - \projH \nabla u_h\right) }{  \nabla w} = \\&= \sum_{E\in \Mh}\scal[E]{\K_E\left(  \nabla u - \projH \nabla u_h\right)}{ \nabla w- \nabla w_I } + \scal[\Omega]{f-f_h}{ w_I}.
\end{align*}
Then, applying \eqref{eq:contVarForm}, it follows
\begin{align}
  \notag
  &\sum_{E\in \Mh} \scal[E]{\K_E \left( \nabla u - \projH \nabla u_h\right) }{ \nabla w} 
  \\
  \notag
  & =  \scal[\Omega]{f}{  w-w_I } - \sum_{E\in \Mh}\scal[E]{\K_E \projH \nabla u_h}{  \nabla w- \nabla w_I}+ \scal[\Omega]{f-f_h}{ w_I}
  \\
  \notag
  & =  \scal[\Omega]{f-f_h}{  w-w_I } +\scal[\Omega]{f_h}{  w-w_I }
  \\
  \notag
  &\quad \quad -\sum_{E\in \Mh}\scal[E]{ \K_E \projH \nabla u_h}{ \nabla w- \nabla w_I }+ \scal[\Omega]{f-f_h}{ w_I}
  \\
  \label{eq:upperBound:firstStep}
  & =  \scal[\Omega]{f_h}{  w-w_I } -\sum_{E\in \Mh}\scal[E]{ \K_E\projH \nabla u_h}{ \nabla w- \nabla w_I }+ \scal[\Omega]{f-f_h}{ w}.
\end{align}
Considering the first two terms above, applying Green's theorem and since
$w-w_I\in\sobh[0]{1}{\Omega}$, we have
\begin{align*}
  &\scal[\Omega]{f_h}{  w-w_I }-\sum_{E\in \Mh}\scal[E]{ \K_E \projH \nabla u_h}{ \nabla w- \nabla w_I }
  \\
  &\qquad= \sum_{E\in \Mh}\scal[E]{f_h+\divt \left( \K_E \, \projH \nabla u_h\right)}{w-w_I}
  \\
  &\qquad\quad - \sum_{e\in \Ehint}\scal[e]{[[\K \, \projH \nabla u_h]]_e}{w-w_I}
  \\ 
  &\qquad= \sum_{E\in \Mh} \scal[E]{r_E}{w-w_I} -
    \sum_{e\in  \Ehint}\scal[e]{j_e}{w-w_I}
  \\
  &\qquad =  I + II 
\end{align*}
By employing \eqref{eq:BoundClementElementK} and with the definition
$C_{\K}\coloneqq \max_{E\in\Mh} \{C_{\K,E}\}$, we get
\begin{align*}
  I& =  \sum_{E\in \Mh}\scal[E]{r_E}{w -w_I}\le \sum_{E\in \Mh}\norm[0,E]{r_E} \norm[0,E]{w-w_I}
  \\&\le C_{I}  \sum_{E\in \Mh}C_{\K, E} \frac{h_E}{\sqrt{\K_E}} \norm[0,E]{r_E} \norm[0,\tilde{\omega}_E]{\sqrt{\K}\nabla w}
  \\  &\le C_{I}  C_{\K}\sum_{E\in \Mh} \frac{h_E}{\sqrt{\K_E}} \norm[0,E]{r_E} \norm[0,\tilde{\omega}_E]{\sqrt{\K}\nabla w}.
\end{align*}
Using the H\"older inequality, we have
\begin{align*}
    I \le C_{I}  C_{\K} \left(\sum_{E\in \Mh} \frac{h^2_E}{\K_E} \norm[0,E]{r_E}^2 \right)^{1/2}\left(\sum_{E\in \Mh}\norm[0,\tilde{\omega}_E]{\sqrt{\K}\nabla w}^2 \right)^{1/2}.
\end{align*}
On the other hand, employing \eqref{eq:BoundClementSqrtEdgeK} and the H\"older inequality, we obtain
\begin{align*}
  II &= -\sum_{e \in \Ehint}  \scal[e]{j_e}{w - w_I}\le \sum_{e \in\Ehint}  \norm[0,e]{j_e} \norm[0,e]{w -w_I}
  \\
     &\le C_{I,1} \sum_{e \in\Ehint}   C_{\K_{\omega_e}} \frac{h_e^{1/2}}{\sqrt{\K_{\omega_e}}}\norm[0,e]{j_e}\norm[0, \tilde{\omega}_e]{\sqrt{\K}\nabla w}
  \\
     &\le C_{I,1}  C^{\Ehint}_{\K}  \sum_{e \in\Ehint}  \frac{h_e^{1/2}}{\sqrt{\K_{\omega_e}}}\norm[0,e]{j_e}\norm[0, \tilde{\omega}_e]{\sqrt{\K}\nabla w}\\
     &\le C_{I,1}  C^{\Ehint}_{\K}  \left(\sum_{e \in\Ehint}
       \frac{h_e}{\K_{\omega_e}}\norm[0,e]{j_e}^2\right)^{1/2}
       \left(\sum_{e \in\Ehint}\norm[0, \tilde{\omega}_e]{\sqrt{\K}\nabla w}^2
       \right)^{1/2},
\end{align*}
where $C^{\Ehint}_{\K} \coloneqq \max_{e\in\Ehint} \left\{ C_{\K_{\omega_e}}\right\}$. We notice that from the mesh quality assumptions and the Definition \ref{def:mesh_variables} we have that
\begin{align*}
  \left(\sum_{E\in \Mh}\norm[0,\tilde{\omega}_E]{\sqrt{\K}\nabla w}^2 \right)^{1/2}
  &\le \sqrt{\Nomegamax}\norm[0,\Omega]{\sqrt{\K}\nabla w},
  \\
  \left(\sum_{e \in\Ehint}\norm[0, \tilde{\omega}_e]{\sqrt{\K}\nabla w}^2
  \right)^{1/2}
  &\le\sqrt{2 \Nomegamax} \norm[0,\Omega]{\sqrt{\K}\nabla w},
\end{align*}
which bring the previous bounds to be
\begin{align}
    I&\le \sqrt{\Nomegamax} C_{I}  C_{\K} \left(\sum_{E\in \Mh} \frac{h^2_E}{\K_E} \norm[0,E]{r_E}^2 \right)^{1/2} \norm[0,\Omega]{\sqrt{\K}\nabla w},\\
    II&\le \sqrt{2 \Nomegamax}  C_{I,1}  C^{\Ehint}_{\K} \left(\sum_{e \in\Ehint}  \frac{h_e}{\K_{\omega_e}}\norm[0,e]{j_e}^2\right)^{1/2} \norm[0,\Omega]{\sqrt{\K}\nabla w} \,.
\end{align}
Going back to \eqref{eq:upperBound:firstStep}, the last term can be estimated
employing Lemma \ref{lemma:projector} and H\"older inequality, as follows:
\begin{align*}
    \sum_{E\in\Mh}\scal[E]{f-f_h}{ w}&= \sum_{E\in\Mh}\scal[E]{f-f_h}{ w - \proj{k}{E}w }  \\
    &\le C_p\sum_{E\in\Mh}\frac{h_E}{\sqrt{\K_E}}\norm[0,E]{f-f_h}\norm[0,E]{\sqrt{\K}\nabla w}\\
    &\le C_p \left(\sum_{E\in\Mh}\frac{h^2_E}{\K_E}\norm[0,E]{f-f_h}^2\right)^{1/2}\left(\sum_{E\in\Mh} \norm[0,E]{\sqrt{\K}\nabla w}^2\right)^{1/2}\\
    & = C_p \left(\sum_{E\in\Mh}\frac{h^2_E}{\K_E}\norm[0,E]{f-f_h}^2\right)^{1/2}\norm[0,\Omega]{\sqrt{\K}\nabla w}.
\end{align*}    
Summing up all the terms we get
\begin{align*}
 & \ennorm[\K,1,\Mh]{u-u_h}=\frac{ \sum_{E\in \Mh} \scal[E]{\K_E \left( \nabla u - \projH \nabla u_h\right) }{ \nabla w}}{\norm[0,\Omega]{\sqrt{\K}\nabla w} }\\
 &  \le \sqrt{\Nomegamax} C_{I}  C_{\K} \left(\sum_{E\in \Mh} \frac{h^2_E}{\K_E} \norm[0,E]{r_E}^2 \right)^{1/2} \\&\quad+  \sqrt{2 \Nomegamax}  C_{I,1}  C^{\Ehint}_{\K} \left(\sum_{e \in\Ehint}  \frac{h_e}{\K_{\omega_e}}\norm[0,e]{j_e}^2\right)^{1/2}+ C_p \left(\sum_{E\in\Mh}\frac{h^2_E}{\K_E}\norm[0,E]{f-f_h}^2\right)^{1/2}.
\end{align*}
Then, the thesis is obtained as follows:
\begin{align*}
  &\ennorm[\K,1,\Mh]{u-u_h}^2\\
  &\qquad\le 3 \,C_u \sum_{E\in \Mh} \left(\frac{h^2_E}{\K_E}
    \norm[0,E]{r_E}^2 +  \sum_{e \in\Eh[E]\cap\Ehint} \frac{h_e}{\K_{\omega_e}}
    \norm[0,e]{j_e}^2+  \frac{h^2_E}{\K_E} \norm[0,E]{f-f_h}^2\right)
  \\
  &\qquad \le 3 \, C_u  \sum_{E\in \Mh}\left( \eta^2_E +\mathcal{F}_E^2\right),
\end{align*}
where $C_u=\left(\max\{\sqrt{\Nomegamax} C_{I}  C_{\K}, \sqrt{2 \Nomegamax}  C_{I,1}  C^{\Ehint}_{\K} ,C_p\}\right)^2$, and $C_U~=~ 3\, C_u$.
\end{proof}
\begin{remark} \label{rmk:quasi-mon} Following
  \cite{Petzoldt2002}, assuming a quasi-monotonicity property  on $\K$, it is possible to define a Cl\'ement-type
  quasi-interpolation operator satisfying error estimates that are independent
  of the local jumps of $\K$. In particular, in the proof of
  \eqref{eq:BoundClementElementK} we get
  \begin{align*}
  &\norm[0, \tilde{\omega}_E]{\nabla v} \le \frac{1}{\sqrt{\K_E}}\norm[0, \tilde{\omega}_E]{\sqrt{\K}\nabla v}, 
\end{align*}
and then $C_{\K,E}=1$. Consequently in \eqref{eq:BoundClementSqrtEdgeK} we have
$C_{\K_{\omega_e}}=1$ and in
Theorem \ref{theorem:upperBound}, the constant $C_U$ does not depend on the
diffusion $\K$.
\end{remark}

\subsection{Lower bound of the error}
In order to prove the lower bound, we use the bubble function
$\psi_E\in \sobh[0]{1}{E}$ for each $E\in\Mh$, as defined in
\cite{berrone_residual_2017,Cangiani2017}. In particular, we build a shape-regular sub-triangulation of $E$ and define $\psi_E$ as the sum of the barycentric bubble functions, which are polynomials on each sub-triangle. As done in \cite{Berrone2006}, for any element $E\in\Mh$,
we define the function
\begin{align}\label{def:w_r_E}
    w_{r,E}(x) \coloneqq \begin{cases}
      \frac{h_E^2}{\K_E}  r_E(x) \psi_E(x)& x \in E \,,
      \\
       0 & x \in \Omega \setminus E \,,
    \end{cases}
\end{align}
where $r_E$ is defined by \eqref{eq:defRE}. Since $r_E$ is a polynomial, using the techniques in
\cite{Verfuerth1996, Berrone2006}, the following results hold true.
\begin{lemma}
  \label{lemma_interior_bubble}
  Let $E\in\Mh$ and $w_{r,E}$ the corresponding function defined above.  The
  following inequalities hold true
  \begin{align*}
    \frac{h_E^2}{\K_E} \norm[0,E]{r_E}^2 \le C_{1,B}\scal[E]{r_E}{w_{r,E}},
  \end{align*}
  and
  \begin{align*}
    h_E \norm[0,E]{\nabla w_{r,E}} \le C_{2,B} \frac{h^2_E}{\K_E}\norm[0,E]{r_E},
  \end{align*}
  where $C_{1,B}$, $C_{2,B}$ are constants independent of $h_E$, but depending on the mesh regularity..
\end{lemma}
\begin{remark}
  \label{rmk:pol_degree}
  From the definition of the space \eqref{eq:firstDecompositionPkell}, it is immediate to check that
  ${\divt\projH}\nabla u \in\Poly{k-2}{E}$. Thus, since $f_h\in\Poly{k}{E}$, $r_E\in\Poly{k}{E}$, which implies that the
  constants $C_{1,B}$ and $C_{2,B}$ in Lemma \ref{lemma_interior_bubble} do not depend on $\ell_E$.
\end{remark}

Similarly, for any given $e\in\Ehint$, let $\psi_e$ be the bubble function
relative to $e$, as defined in \cite{berrone_residual_2017,Cangiani2017}. In particular, we consider the sub-triangles sharing $e$ of the elements in $\omega_e$ and $\psi_e$ is defined as the sum of the barycentric bubble functions relative to these sub-triangles. 
Moreover, following \cite{Verfuerth1996} we extend $j_e$, defined by
\eqref{eq:defJE}, through a constant prolongation in the normal direction with
respect to $e$. Let $\mathcal{C}j_e$ be such function. We define
\begin{align}\label{def:_j_E}
  w_{j,e}(x) \coloneqq
  \begin{cases}
    \frac{h_e}{\K_{\omega_e}} (\mathcal{C}j_e)(x)\psi_e(x)& x \in \omega_e,
    \\
    0 & x \in \Omega \setminus \omega_e\,.
  \end{cases} 
\end{align}
The following results can be proved using the techniques in \cite{Verfuerth1996, Berrone2006}.
\begin{lemma}
  \label{lemma_edge_bubble}
  Let $E\in\Mh$, $e\in \Eh[E]$ and $w_{j,E}$ as defined in
  \eqref{def:_j_E}. Then,
  \begin{align*}
    \frac{h_e}{\K_{\omega_e}} \norm[0,e]{j_e}^2\le C_{1,b}\scal[e]{j_e}{ w_{j,e}},
  \end{align*}
  and
  \begin{align*}
    h_e^{1/2} \norm[0,E]{\nabla w_{j,e}} \le C_{2,b}  \frac{h_e}{\K_{\omega_e}} \norm[0,e]{j_e},
  \end{align*}
  where $C_{1,b}$, $C_{2,b}$ are constants independent of $h_E$, but depending on the mesh regularity.
\end{lemma}
\begin{theorem}[Local lower bound] \label{theorem:lowerBound}

  Let $E\in\Mh$ be given, $u$ the solution to the continuous problem \eqref{eq:contVarForm} and
  $u_h \in\V{k}$ the solution to the discrete problem \eqref{eq:higher_order_discrVarForm}.  Then
  there exists a constant $C_L>0$ independent of $h_E$, but depending on the mesh regularity, such that
  \begin{align*}
    \eta^2_E\le C_L  \left( \ennorm[\K,1,\omega_E]{u-u_ h}^2 +
    \sum_{E\in \omega_E} \mathcal{F}^2_E \right),
  \end{align*}
  where $\ennorm[\K,1,\omega_E]{u-u_h }$ is defined by \eqref{def:error} and $\eta_E$ and
  $\mathcal{F}_E$ in Theorem \ref{theorem:upperBound}.
\end{theorem}
\begin{proof}
  For any $w\in\sobho{\Omega}{}$, exploiting problems \eqref{eq:contVarForm} and
  \eqref{eq:higher_order_discrVarForm} and the definitions of $r_E$ and $j_e$ in
  \eqref{eq:defRE}-\eqref{eq:defJE}, we have
\begin{align}
  \notag
  &\sum_{E\in \Mh} \scal[E]{\K_E \left( \nabla u - \projH \nabla u_h\right) }{ \nabla w}=
  \\
  \notag
  &= \sum_{E\in \Mh} \scal[E]{f}{ w} -\sum_{E\in \Mh} \scal[E]{\K_E  \projH \nabla u_h }{\nabla w}
  \\
  \notag
  &= \sum_{E\in \Mh} \scal[E]{f_h}{ w} -\sum_{E\in \Mh} \scal[E]{\K_E  \projH \nabla u_h }{ \nabla w} +\sum_{E\in \Mh} \scal[E]{f-f_h}{ w}
  \\
  \notag
  &= \sum_{E\in \Mh} \scal[E]{f_h}{ w} +\sum_{E\in \Mh} \scal[E]{\divt\left (\K_E  \projH \nabla u_h \right)}{ w}
  \\
  \notag
  & \qquad-\sum_{e\in \Ehint}
    \scal[e]{[[\K \, \projHOmega \nabla u_h ]]_e}{ w} +\sum_{E\in \Mh} \scal[E]{f-f_h}{ w}
  \\
  \label{eq:lowerbound_firstEquation}
  & = \sum_{E\in \Mh} \scal[E]{r_E}{ w} -\sum_{e\in \Ehint} \scal[e]{j_e}{ w} +\sum_{E\in \Mh} \scal[E]{f-f_h}{ w} \,.
\end{align}
First, taking $w= w_{r,E}$ as defined in \eqref{def:w_r_E}, since $\supp(w_{r,E}) \subseteq E$ we
get
 \begin{equation*}
\scal[E]{\K_E \left( \nabla u - \projH \nabla u_h\right) }{ \nabla w_{r,E}}= \scal[E]{f-f_h}{ w_{r,E}} + \scal[E]{r_E}{ w_{r,E} }.  
 \end{equation*}
Using the Cauchy-Schwarz inequality and Lemma \ref{lemma_interior_bubble}, we have that
\begin{align*}
  & \frac{h^2_E}{\K_E}\norm[0,E]{r_E}^2 \le  C_{1,B}\scal[E]{r_E}{ w_{r,E}}
  \\
  &= C_{1,B} \scal[E]{\K_E \left( \nabla u - \projH \nabla u_h\right) }{ \nabla w_{r,E}}+ C_{1,B}\scal[E]{f_h-f}{ w_{r,E}}
  \\
  &\le C_{1,B} \scal[E]{\K_E \left( \nabla u - \projH \nabla u_h\right) }{ \nabla w_{r,E}}+C_{1,B}\norm[0,E]{f-f_h}\norm[0,E]{ w_{r,E}}
  \\
  &\le C_{1,B} \scal[E]{\K_E \left( \nabla u - \projH \nabla u_h\right) }{ \nabla w_{r,E}}
  \\&\qquad+C_{1,B} C_{p} \frac{h_E}{\sqrt{\K_E}}\norm[0,E]{f-f_h}\norm[0,E]{\sqrt{\K}\nabla w_{r,E}}
  \\
  &= C_{1,B} \scal[E]{\K_E \left( \nabla u - \projH \nabla u_h\right) }{\nabla w_{r,E}}
    \\
    &\qquad+C_{1,B} C_{p}\left( \frac{h^2_E}{\K_E}\norm[0,E]{f-f_h}^2\right)^{1/2}\norm[0,E]{\sqrt{\K}\nabla w_{r,E}} \,,
\end{align*}
where $C_{p}$ is a Poincar\'e constant.
Then,
\begin{equation}
   \label{eq:lowerBound:residual-step}
 \frac{\frac{h^2_E}{\K_E}\norm[0,E]{r_E}^2}{\norm[0,E]{\sqrt{\K}\nabla w_{r,E}}}  \le C_{1,B} \ennorm[\K,1,E]{u-u_h} + C_{1,B} C_{p}\left(\frac{h^2_E}{\K_E}\norm[0,E]{f-f_h}^2\right)^{1/2}.
\end{equation}
We notice that, from Lemma \ref{lemma_interior_bubble}, one has
\begin{equation*}
  \norm[0,E]{\sqrt{\K}\nabla w_{r,E}} =\sqrt{\K_E}
    \norm[0,E]{\nabla w_{r,E}}
    \le C_{2,B}  \frac{h_E}{\sqrt{\K_E}} \norm[0,E]{ r_E}.
\end{equation*}
Then, from \eqref{eq:lowerBound:residual-step} we get
\begin{equation*}
   \frac{h_E}{\sqrt{\K_E}} \norm[0,E]{ r_E} \le C_{1,B} C_{2,B} \ennorm[\K,1,E]{u-u_h} + C_{1,B} C_{2,B}C_{p}\left( \frac{h^2_E}{\K_E}\norm[0,E]{f-f_h}^2\right)^{1/2}, 
\end{equation*}
and, squaring,
\begin{align}\label{eq:boundResidual_lower}
     \frac{h^2_E}{\K_E} \norm[0,E]{r_E}^2 \le  C_{1,L}\left( \ennorm[\K,1,E]{u-u_h}^2+ \frac{h^2_E}{\K_E}\norm[0,E]{f-f_h}^2\right),
\end{align} 
where $C_{1,L}= 2 (\max\{ C_{1,B} C_{2,B}, C_{1,B} C_{2,B} C_{p}\})^2$.  Going back to
\eqref{eq:lowerbound_firstEquation}, for any edge $e\in\Ehint$ we take $w=w_{j,e}$ as defined in
\eqref{def:_j_E}. Since $\supp(w_{j,e})\subseteq \omega_e$ we
get
\begin{equation*}
\begin{split}
    \sum_{E\in \omega_e} \scal[E]{\K_E\left( \nabla u - \projH \nabla u_h\right) }{ \nabla w_{j,e}}
  &= \sum_{E\in \omega_e} \scal[E]{r_E}{ w_{j,e}} - \scal[e]{j_e}{w_{j,e}}
  \\
  &\qquad+\sum_{E\in \omega_e} \scal[E]{f-f_h}{w_{j,e}}.
\end{split}
\end{equation*}
Using Lemma \ref{lemma_edge_bubble} and the H\"older and Poincar\'e inequalities, we get
\begin{align*}
  \frac{h_e}{\K_{\omega_e}}\norm[0,e]{j_e}^2
  &\le C_{1,b}\scal[e]{j_e}{ w_{j,e}} = C_{1,b}\sum_{E\in\omega_e} \left(\scal[E]{r_E}{w_{j,e}} +\scal[E]{f-f_h}{w_{j,e}}\right)
  \\
  &\qquad- C_{1,b}\sum_{E\in\omega_e}\scal[E]{\K_E \left( \nabla u - \projH \nabla u_h\right)} { \nabla w_{j,e}}
  \\
  &\le  C_{1,b} C_{p}\sum_{E\in\omega_e}
    \frac{h_E}{\sqrt{\K_E}}\norm[0,E]{r_E}\norm[0,E]{\sqrt{\K} \nabla w_{j,e}}
  \\
  &\quad + C_{1,b} C_{p}\sum_{E\in\omega_e}
    \frac{h_E}{\sqrt{\K_E}}\norm[0,E]{f-f_h} \norm[0,E]{ \sqrt{\K}\nabla w_{j,e}}
  \\
  &\quad - C_{1,b} \sum_{E\in\omega_e}
    \scal[E]{\K_E \left( \nabla u - \projH \nabla u_h\right)}{\nabla w_{j,e}}
  \\
  &\le  C_{1,b} C_{p}\left(\sum_{E\in\omega_e} \frac{h^2_E}{\K_E}\norm[0,E]{r_E}^2 \right)^{1/2}\left(\sum_{E\in\omega_e}\norm[0,E]{\sqrt{\K}\nabla w_{j,e}}^2\right)^{1/2} 
  \\
  &\quad + C_{1,b} C_{p}\left(\sum_{E\in\omega_e}\frac{h^2_E}{\K_E}\norm[0,E]{f-f_h}^2 \right)^{1/2}\left(\sum_{E\in\omega_e}\norm[0,E]{\sqrt{\K}\nabla w_{j,e}}^2 \right)^{1/2}
  \\
  &\quad - C_{1,b} \sum_{E\in\omega_e} \scal[E]{\K_E \left( \nabla u - \projH \nabla u_h\right)} { \nabla w_{j,e}}.
\end{align*}
Then, employing the definition of $\ennorm[\K,1,\omega_e]{u-u_h}$ in \eqref{def:error},
\begin{equation}
  \label{eq:lowerBound:jump-estimate}
  \begin{split}
    &\frac{\frac{h_e}{\K_{\omega_e}}
      \norm[0,e]{j_e}^2}{\left(\sum_{E\in\omega_e}
      \norm[0,E]{\sqrt{\K}\nabla w_{j,e}}^2\right)^{1/2}}
    \\
    &\le  C_{1,b} C_{p} \left(\sum_{E\in\omega_e} \frac{h^2_E}{\K_E} \norm[0,E]{r_E}^2 \right)^{1/2}
      + C_{1,b} C_{p}\left(\sum_{E\in\omega_e}\frac{h^2_E}{\K_E}\norm[0,E]{f-f_h}^2 \right)^{1/2}
    \\
    & \quad + C_{1,b} \frac{\sum_{E\in\omega_e}
      \scal[E]{\K_E \left( \nabla u - \projH \nabla u_h\right)}{ \nabla \left(-w_{j,e}\right)}}{\left(\sum_{E\in\omega_e}\norm[0,E]{\sqrt{\K}\nabla \left(-w_{j,e}\right)}^2\right)^{1/2}}
    \\
    &\le  C_{1,b} C_{p} \left(\sum_{E\in\omega_e} \frac{h^2_E}{\K_E} \norm[0,E]{r_E}^2 \right)^{1/2}
      + C_{1,b} C_{p}\left(\sum_{E\in\omega_e}\frac{h^2_E}{\K_E}\norm[0,E]{f-f_h}^2 \right)^{1/2}
    \\
    & \quad + C_{1,b} \ennorm[\K,1,\omega_e]{u-u_h} \,.
  \end{split}
\end{equation}
We employ again Lemma \ref{lemma_edge_bubble}, obtaining
\begin{align*}
  \left(\sum_{E\in\omega_e}\norm[0,E]{\sqrt{\K}\nabla w_{j,e} }^2\right)^{1/2}
  &\le C_{2,b} \left(\sum_{E\in\omega_e} \K_E\frac{h_e}{\K^2_{\omega_e}}
  \norm[0,e]{ j_e}^2\right)^{1/2}\\
  &\qquad\le C_{2,b} \left(\frac{h_e}{\K_{\omega_e}}\norm[0,e]{ j_e}^2\right)^{1/2} \,.
\end{align*}
Thus, subtituting in \eqref{eq:lowerBound:jump-estimate}, we get
\begin{equation*}%\label{eq:boundJump_lower1}
\begin{split}
  \left(\frac{h_e}{\K_{\omega_e}}\norm[0,e]{ j_e}^2\right)^{1/2}
  \le C_{1,b}  C_{2,b} C_{p}\left(\sum_{E\in\omega_e} \frac{h^2_E}{\K_E} \norm[0,E]{r_E}^2 \right)^{1/2}
  \\
  + C_{1,b} C_{2,b} C_{p} \left(\sum_{E\in\omega_e}\frac{h^2_E}{\K_E}\norm[0,E]{f-f_h}^2 \right)^{1/2}+  C_{1,b} C_{2,b} \ennorm[\K,1,\omega_e]{u-u_h}.
\end{split}
\end{equation*}
and squaring,
\begin{align}\label{eq:boundJump_lower}
    &\frac{h_e}{\K_{\omega_e}}\norm[0,e]{ j_e}^2\le C_{2, L}\left(\sum_{E\in\omega_e} \frac{h^2_E}{\K_E} \norm[0,E]{r_E}^2 +\sum_{E\in\omega_e}\frac{h^2_E}{\K_E}\norm[0,E]{f-f_h}^2 + \ennorm[\K,1,\omega_e]{u-u_h}^2\right).
\end{align}
where $C_{2,L}= 3 \max\{ C_{1,b} C_{2,b}, C_{1,b} C_{2,b} C_{p} \}^2 $.
Finally, we apply \eqref{eq:boundResidual_lower} in \eqref{eq:boundJump_lower}, getting
\begin{align}\label{eq:boundJump_lower_FINAL}
    \frac{h_e}{\K_{\omega_e}}\norm[0,e]{j_e}^2
    &\le C_{3, L} \left( \sum_{E\in \omega_e}\frac{h^2_{E}}{\K_{E}}\norm[0,E]{f-f_h}^2+\ennorm[\K,1,\omega_e]{u-u_h}^2
    \right),
\end{align}
where $C_{3,L}= (C_{1,L}+1)C_{2,L}$.
We conclude the proof by summing \eqref{eq:boundResidual_lower} and \eqref{eq:boundJump_lower_FINAL}. 
\end{proof}
\begin{corollary}[Global lower bound]
  Let $u$ the solution to the continuous problem \eqref{eq:contVarForm} and $u_h \in\V{k}$ the
  solution to the discrete problem \eqref{eq:higher_order_discrVarForm}.\\ Let
  $ \ennorm[\K,1,\Mh]{u-u_h}$ be defined as in \eqref{def:error}, it holds
  \begin{align*}
    \sum_{E\in \Mh}\eta^2_E\le C_L  (\NVmax+1) \left( \ennorm[\K,1,\Mh]{u-u_ h}^2 +\sum_{E\in \Mh} \mathcal{F}^2_E \right),
  \end{align*}  
  where $\eta_E, \,\mathcal{F}_E$ are defined in Theorem \ref{theorem:upperBound}, $C_L$ is defined
  in Theorem \ref{theorem:lowerBound}, and $\NVmax$ as in Assumption \ref{hp:mesh_assumption_VEM}.
\end{corollary}

\section{Numerical results} \label{sec:NumericalResults}

In this section, we present some numerical tests to show the equivalence between error and error estimator in the different problems proposed in \cite{berrone_residual_2017} and on a standard test
on an L-shaped domain. Firstly, we approximate the error defined in \eqref{def:error} as
\begin{align*}
  \ennorm[\K,1,\Mh]{u-u_h} \simeq
  \left(
  \sum_{E\in\Mh}\norm[0,E]{\sqrt{\K
  }  \left( \nabla u - \projH \nabla u_h\right) }^2
  \right)^{\frac12},
\end{align*}
and define the \emph{effectivity index} as the ratio between the estimator and the error, i.e.
\begin{align*}
  \epsilon \coloneqq \left(\frac{\sum_{E\in\Mh}\eta^2_{E}}
  {
  \sum_{E\in\Mh}\norm[0,E]{\sqrt{\K
  }  \left( \nabla u - \projH \nabla u_h\right) }^2
  }\right)^{\frac12}.
\end{align*}
\subsection{Test 1}
\begin{figure}[!ht]  
\centering 
\begin{subfigure}[b]{.3\linewidth} \centering
    \includegraphics[width=1.\linewidth]{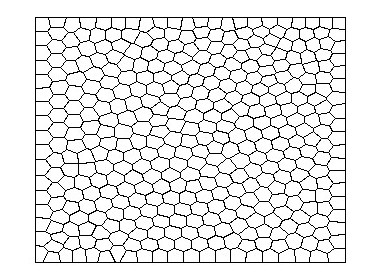}
    \caption{\meshtag{Polymesher}}
    \label{subfig:Polymesher_mesh}
  \end{subfigure} \quad
\begin{subfigure}[b]{.3\linewidth} \centering
    \includegraphics[width=1.\linewidth]{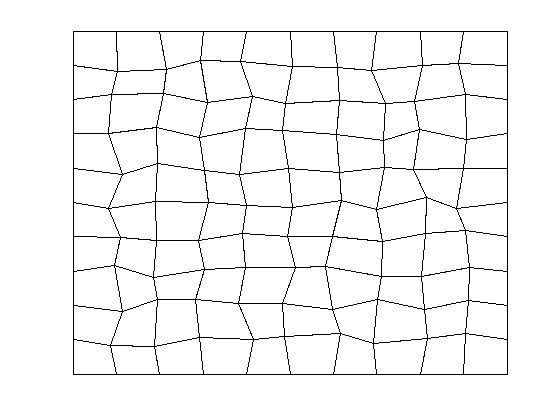}
    \caption{\meshtag{Distorted cartesian}}
    \label{subfig:DistortedCartesian_mesh}
  \end{subfigure}  \quad
  \begin{subfigure}[b]{.3\linewidth} \centering
    \includegraphics[width=1.\linewidth]{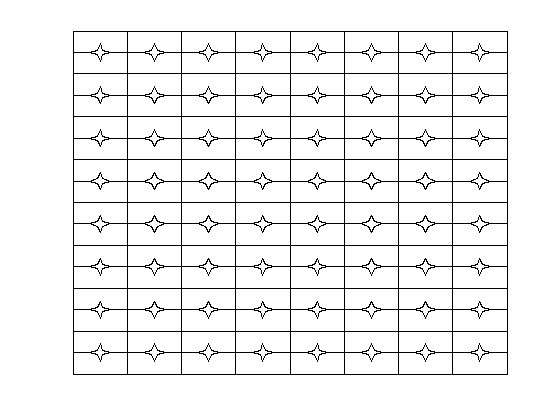}
    \caption{\meshtag{Star concave}}
    \label{subfig:starconcave_mesh}
  \end{subfigure}
\caption{Meshes used for convergence tests.}
  \label{fig:Meshes}
  \end{figure}

\begin{figure}[!ht]  \centering

\begin{subfigure}[b]{.3\linewidth} \centering
    \includegraphics[width=1.\linewidth]{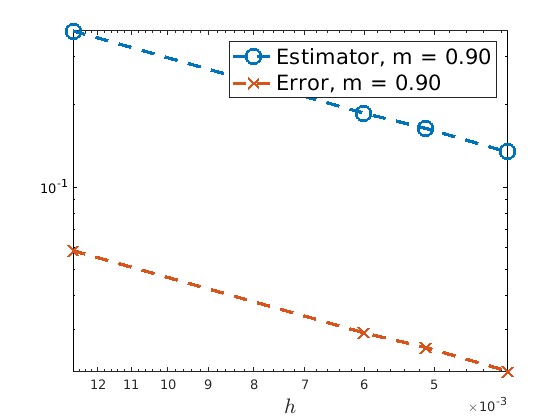}
    \caption{\meshtag{Polymesher}, $k=1$.}
    \label{subfig:Polymesher_k1}
  \end{subfigure}\quad
 \begin{subfigure}[b]{.3\linewidth} \centering
    \includegraphics[width=1.\linewidth]{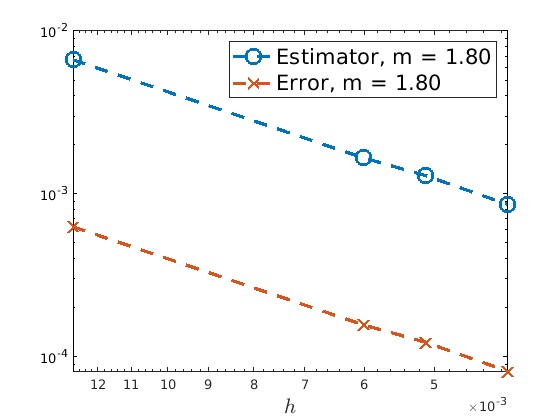}
    \caption{\meshtag{Polymesher}, $k=2$.}
    \label{subfig:Polymesher_k2}
  \end{subfigure}\quad
  \begin{subfigure}[b]{.3\linewidth} \centering
    \includegraphics[width=1.\linewidth]{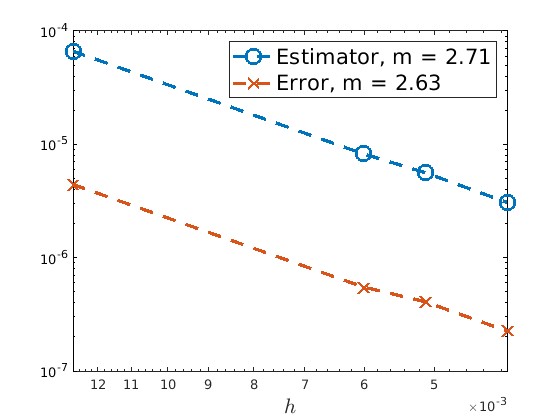}
    \caption{\meshtag{Polymesher}, $k=3$.}
    \label{subfig:Polymesher_k3}
  \end{subfigure}
  \\
  \begin{subfigure}[b]{.3\linewidth} \centering
    \includegraphics[width=1.\linewidth]{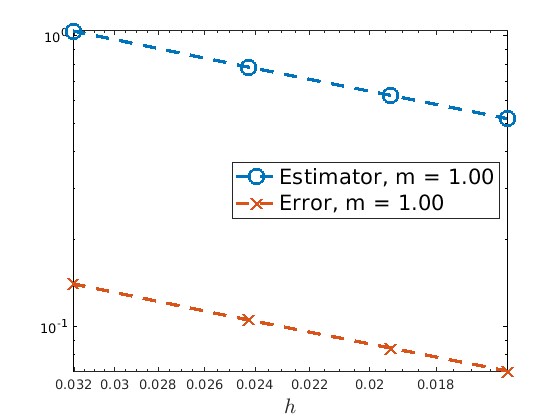}
    \caption{\meshtag{Distorted cartesian}, $k=1$.}
    \label{subfig:DistortedCartesian_k1}
  \end{subfigure}\quad
    \begin{subfigure}[b]{.3\linewidth} \centering
    \includegraphics[width=1.\linewidth]{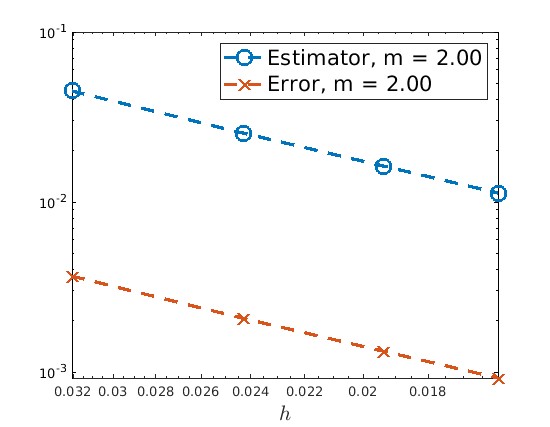}
    \caption{\meshtag{Distorted cartesian}, $k=2$.}
    \label{subfig:DistortedCartesian_k2}
  \end{subfigure}\quad
     \begin{subfigure}[b]{.3\linewidth} \centering
    \includegraphics[width=1.\linewidth]{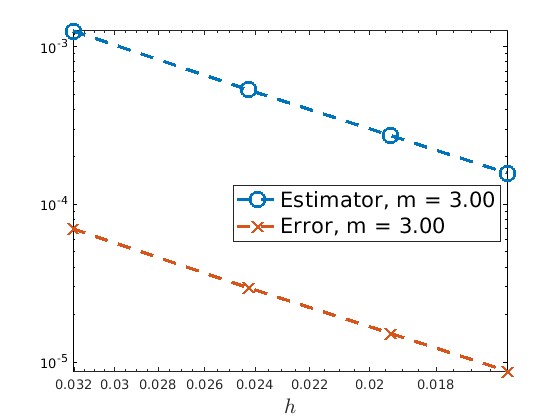}
    \caption{\meshtag{Distorted cartesian}, $k=3$.}
    \label{subfig:DistortedCartesian_k3}
  \end{subfigure}\\
    \begin{subfigure}[b]{.3\linewidth} \centering
    \includegraphics[width=1.\linewidth]{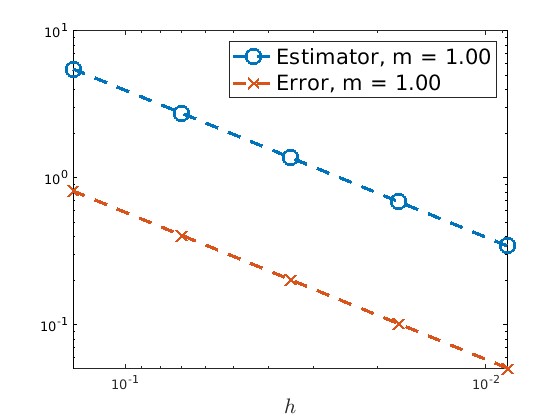}
    \caption{\meshtag{Star concave}, $k=1$.}
    \label{subfig:StarConcave_k1}
  \end{subfigure}\quad
    \begin{subfigure}[b]{.3\linewidth} \centering
    \includegraphics[width=1.\linewidth]{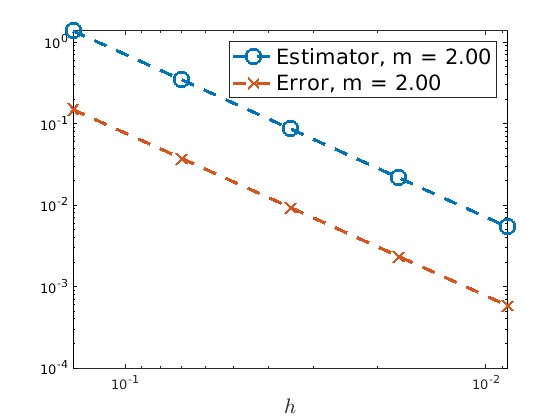}
    \caption{\meshtag{Star concave}, $k=2$.}
    \label{subfig:StarConcave_k2}
  \end{subfigure}\quad
     \begin{subfigure}[b]{.3\linewidth} \centering
    \includegraphics[width=1.\linewidth]{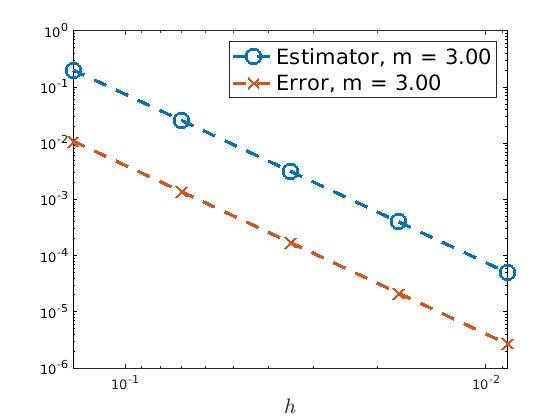}
    \caption{\meshtag{Star concave}, $k=3$.}
    \label{subfig:StarConcave_k3}
  \end{subfigure}
  \caption{Test 1: Convergence plots. In the legends, $m$ is the average convergence rate.}
  \label{Test1_convergence}
  \end{figure}

We consider the Poisson problem with $\K = 1$, setting the loading term $f$ in such a way that the solution to the problem is $u(x, y) = \sin(2 \pi x) \sin(2\pi y)$. We tested it on meshes made up of convex polygons labeled \meshtag{Polymesher} \cite{polymesher}, distorted cartesian meshes (\meshtag{Distorted cartesian}), and concave meshes labeled \meshtag{Star concave}, as represented in Figure \ref{fig:Meshes}. Figure \ref{Test1_convergence} shows the behavior of the estimator and the errors in the cases $k=1,2,3$.
Tables \ref{tab:Polymesher_Test1}, \ref{tab:distorted_Test1} and \ref{tab:StarConcave_Test1} show that the effectivity indices are independent of the meshsize and display a weak dependence on the type of polygons used.
  
\begin{table}[!ht]
    \centering
\begin{tabular}{c c |c c |c c}
\multicolumn{2}{c |}{$k$ = 1} &\multicolumn{2}{c |}{$k$ = 2} &\multicolumn{2}{c}{$k$ = 3} \\ 
\hline 
$h$ & $\epsilon$ &$h$ & $\epsilon$ &$h$ & $\epsilon$ \\ 
0.012807&6.3810 &0.012807&10.6155&0.012807&15.0051 \\ 
0.0060112&6.3832&0.0060112&10.6057&0.0060112&14.9574 \\ 
0.005119&6.3515&0.005119&10.5739&0.005119&13.7805 \\ 
0.0041359&6.3506&0.0041359&10.5712&0.0041359&13.6350 \\ 
\hline 
\end{tabular}
\caption{Test 1, \meshtag{Polymesher}. Effectivity indices.}
    \label{tab:Polymesher_Test1}
\end{table}

\begin{table}[!ht]
    \centering
\begin{tabular}{c c |c c |c c}
\multicolumn{2}{c |}{$k$ = 1} &\multicolumn{2}{c |}{$k$ = 2} &\multicolumn{2}{c}{$k$ = 3} \\ 
\hline 
$h$ & $\epsilon$ &$h$ & $\epsilon$ &$h$ & $\epsilon$ \\ 
0.032068&7.4085&0.032068&12.2719&0.032068&17.9411 \\ 
0.024269&7.4092&0.024269&12.2766&0.024269&17.9522 \\ 
0.019351&7.3965&0.019351&12.2720&0.019351&17.9547 \\ 
0.016079&7.3853&0.016079&12.2419&0.016079&17.9495 \\ 
\hline 
\end{tabular}
\caption{Test 1, \meshtag{Distorted cartesian}. Effectivity indices.}
    \label{tab:distorted_Test1}
\end{table}

\begin{table}[!ht]
    \centering
\begin{tabular}{c c |c c |c c}
\multicolumn{2}{c |}{$k$ = 1} &\multicolumn{2}{c |}{$k$ = 2} &\multicolumn{2}{c}{$k$ = 3} \\ 
\hline 
$h$ & $\epsilon$ &$h$ & $\epsilon$ &$h$ & $\epsilon$ \\ 
0.13975&6.7571&0.13975&9.2299&0.13975&18.4068 \\ 
0.069877&6.7709&0.069877&9.2319&0.069877&18.5379 \\ 
0.034939&6.7736&0.034939&9.2350&0.034939&18.5714 \\ 
0.017469&6.7742&0.017469&9.2360&0.017469&18.5798 \\ 
0.0087346&6.7743&0.0087346&9.2363&0.0087346&18.5720 \\ 
\hline 
\end{tabular} 
\caption{Test 1, \meshtag{Star concave}. Effectivity indices.}
    \label{tab:StarConcave_Test1}
\end{table}

\subsection{Test 2}
We test the estimator in the presence of diffusion jumps. 
Let $\Omega=[0,1]^2$ and the diffusion tensor $\K(x,y) = \gamma_i(x,y) I$, $i=1,2$, and
\begin{align*}
\gamma_1(x,y) \coloneqq
\begin{cases}10,
\;\text{in }\Omega_1 :=  [0, 0.5] \times[0, 1],\\
1, \;\text{in } \Omega_2 :=(0.5, 1] \times[0, 1],
\end{cases} \gamma_2(x,y) \coloneqq
\begin{cases}10^{-3},
\;\text{in } \Omega_1,\\
1, \;\text{in } \Omega_2.
\end{cases}
\end{align*}
The loading term is chosen so that the solution $u_i (x, y) = \xi_i (x)Y (y)$, where
\begin{align}\label{def:xi}
   \xi_i(x) &\coloneqq \begin{cases}
   -\frac{1}{\gamma_i|_{\Omega_1}}\left(\frac{x^2}{2} +c_i x\right)
\quad\text{if  } x\in[0, 0.5],\\
   -\frac{1}{\gamma_i|_{\Omega_2}}\left(\frac{x^2}{2} +c_i x -c_i - \frac{1}{2}\right)
\quad\text{if  } x\in( 0.5,1],
\end{cases}
\end{align}
\begin{align}\label{def:Yy}
    Y(y)&\coloneqq y (1-y) \left(y-\frac{1}{2}\right)^2,
\end{align}
and $c_i\coloneqq - \frac{3 \gamma_i|_{\Omega_1}+ \gamma_{\Omega_2}}{4\left(\gamma_i|_{\Omega_1}+ \gamma_i|_{\Omega_2}\right)}$. 

\begin{figure}[!ht]
  \begin{subfigure}[b]{.3\linewidth} \centering
    \includegraphics[width=1.\linewidth]{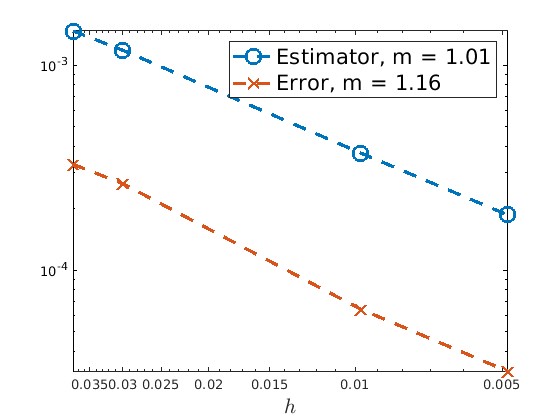}
    \caption{\meshtag{Distorted cartesian}, $\gamma_1$, $k=1$.}
    \label{subfig:DistortedCartesian_k1_Test2}
  \end{subfigure}\quad
    \begin{subfigure}[b]{.3\linewidth} \centering
    \includegraphics[width=1.\linewidth]{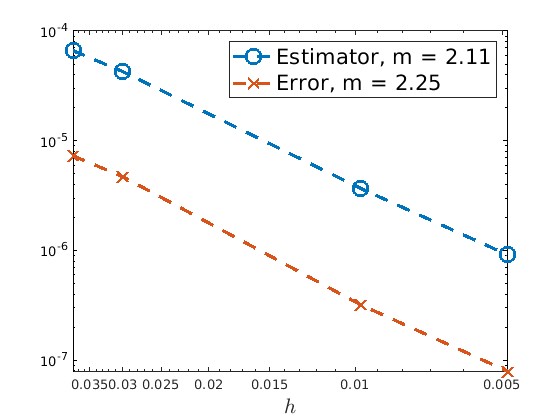}
    \caption{\meshtag{Distorted cartesian}, $\gamma_1$, $k=2$.}
    \label{subfig:DistortedCartesian_k2_Test2}
  \end{subfigure}\quad
     \begin{subfigure}[b]{.3\linewidth} \centering
    \includegraphics[width=1.\linewidth]{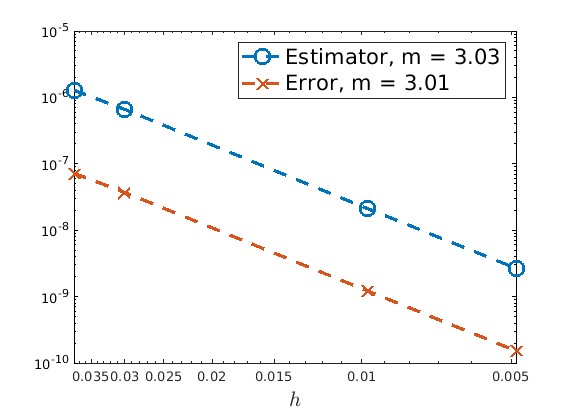}
    \caption{\meshtag{Distorted cartesian},  $\gamma_1$, $k=3$.}
    \label{subfig:DistortedCartesian_k3_Test2}
  \end{subfigure}
   \begin{subfigure}[b]{.3\linewidth} \centering
    \includegraphics[width=1.\linewidth]{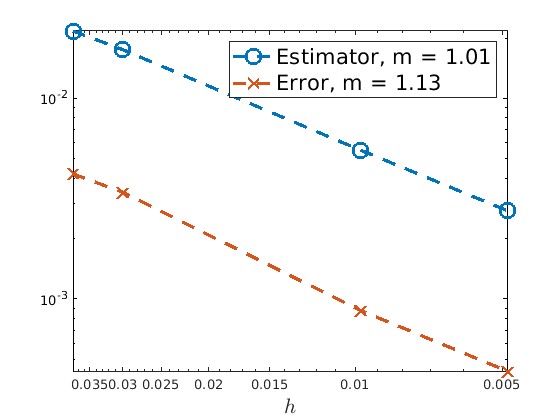}
    \caption{\meshtag{Distorted cartesian}, $\gamma_2$, $k=1$.}
    \label{subfig:DistortedCartesian_k1_Test2bis}
  \end{subfigure}\quad
    \begin{subfigure}[b]{.3\linewidth} \centering
    \includegraphics[width=1.\linewidth]{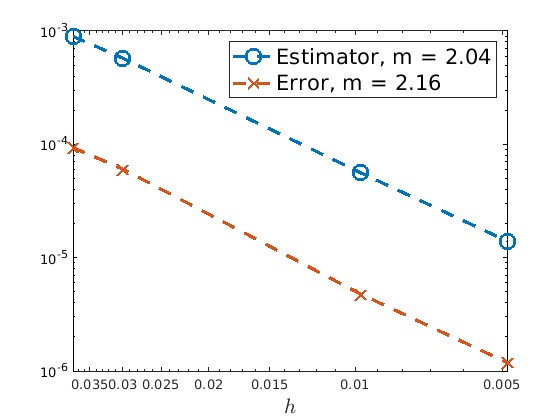}
    \caption{\meshtag{Distorted cartesian}, $\gamma_2$, $k=2$.}
    \label{subfig:DistortedCartesian_k2_Test2bis}
  \end{subfigure}\quad
     \begin{subfigure}[b]{.3\linewidth} \centering
    \includegraphics[width=1.\linewidth]{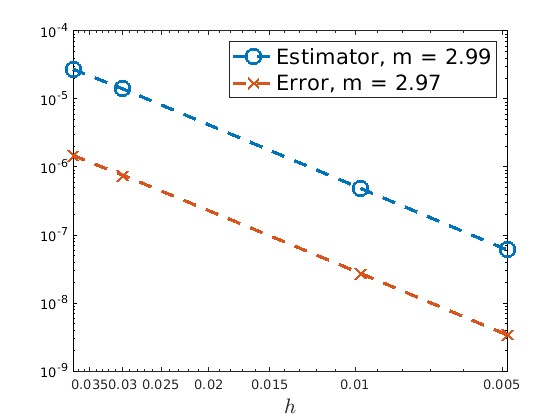}
    \caption{\meshtag{Distorted cartesian},  $\gamma_2$, $k=3$.}
    \label{subfig:DistortedCartesian_k3_Test2bis}
  \end{subfigure}
  \caption{Test 2: Convergence plots. In the legends, $m$ is the average convergence rate.}
  \label{Test2_convergence}
  \end{figure}

The convergence plots displayed in Figure \ref{Test2_convergence} confirm the optimal convergence rates. The effectivity indices remain stable with respect to the meshsize, as shown in Tables \ref{tab:distorted_Test2_gamma1} and \ref{tab:distorted_Test2_gamma2}. The jump of the diffusivity coefficient has a minimal impact on the indices $\epsilon$, demonstrating strong robustness of the estimator.

\begin{table}[!ht]
    \centering
\begin{tabular}{c c |c c |c c}
\multicolumn{2}{c |}{$k$ = 1} &\multicolumn{2}{c |}{$k$ = 2} &\multicolumn{2}{c}{$k$ = 3} \\ 
\hline 
$h$ & $\epsilon$ &$h$ & $\epsilon$ &$h$ & $\epsilon$ \\ 
0.037919&4.4720&0.037919&9.0937&0.037919&17.8213 \\ 
0.029971&4.4795&0.029971&9.0737&0.029971&17.7688 \\ 
0.009748&5.7849&0.009748&11.5067&0.009748&17.4022 \\ 
0.0048901&5.7966&0.0048901&11.5175&0.0048901&17.4179 \\ 
\hline 
\end{tabular}
\caption{Test 2, case $\gamma_1$. \meshtag{Distorted cartesian}. Effectivity indices. }
    \label{tab:distorted_Test2_gamma1}
\end{table}

\begin{table}[!ht]
    \centering
\begin{tabular}{c c |c c |c c}
\multicolumn{2}{c |}{$k$ = 1} &\multicolumn{2}{c |}{$k$ = 2} &\multicolumn{2}{c}{$k$ = 3} \\ 
\hline 
$h$ & $\epsilon$ &$h$ & $\epsilon$ &$h$ & $\epsilon$ \\ 
0.037919&5.1599&0.037919&9.5695&0.037919&18.5862 \\ 
0.029971&5.1422&0.029971&9.5063&0.029971&18.5499 \\ 
0.009748&6.3043&0.009748&11.7771&0.009748&17.7856 \\ 
0.0048901&6.3107&0.0048901&11.7849&0.0048901&17.8039 \\ 
\hline 
\end{tabular}
\caption{Test 2, case $\gamma_2$. \meshtag{Distorted cartesian}. Effectivity indices. }
    \label{tab:distorted_Test2_gamma2}
\end{table}
\subsection{Test 3}
\begin{figure}[!ht]
  \begin{subfigure}[b]{.3\linewidth} \centering
    \includegraphics[width=1.\linewidth]{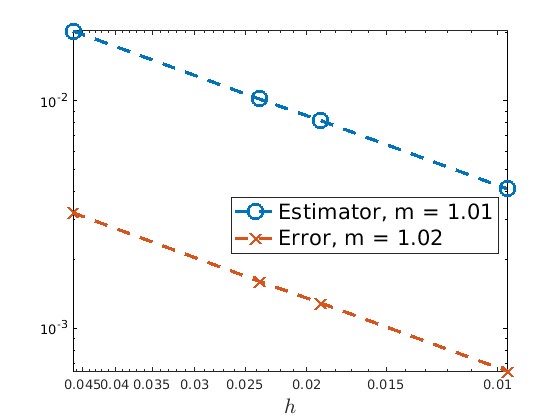}
    \caption{\meshtag{Distorted cartesian}, $\gamma_3$, $k=1$.}
    \label{subfig:DistortedCartesian_k1_Test3}
  \end{subfigure}\quad
    \begin{subfigure}[b]{.3\linewidth} \centering
    \includegraphics[width=1.\linewidth]{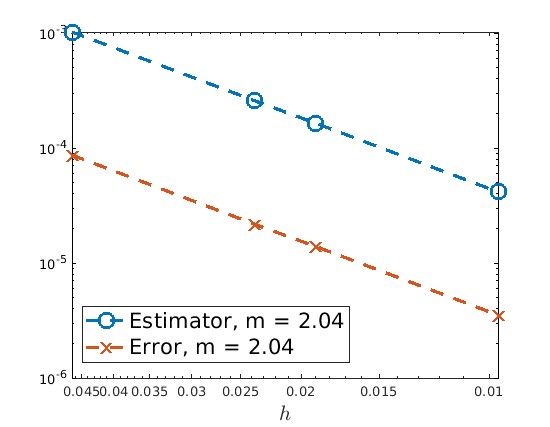}
    \caption{\meshtag{Distorted cartesian}, $\gamma_3$, $k=2$.}
    \label{subfig:DistortedCartesian_k2_Test3}
  \end{subfigure}\quad
     \begin{subfigure}[b]{.3\linewidth} \centering
    \includegraphics[width=1.\linewidth]{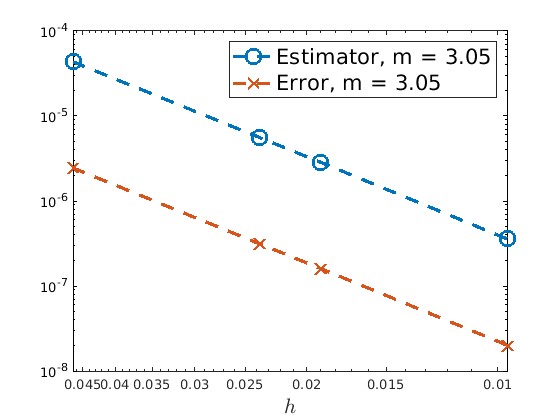}
    \caption{\meshtag{Distorted cartesian},  $\gamma_3$, $k=3$.}
    \label{subfig:DistortedCartesian_k3_Test3}
  \end{subfigure}
   \begin{subfigure}[b]{.3\linewidth} \centering
    \includegraphics[width=1.\linewidth]{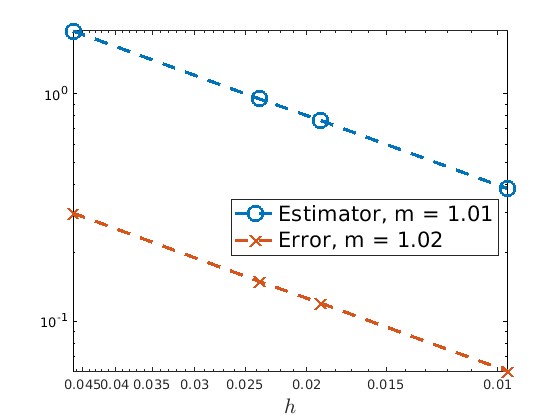}
    \caption{\meshtag{Distorted cartesian}, $\gamma_4$, $k=1$.}
    \label{subfig:DistortedCartesian_k1_Test3bis}
  \end{subfigure}\quad
    \begin{subfigure}[b]{.3\linewidth} \centering
    \includegraphics[width=1.\linewidth]{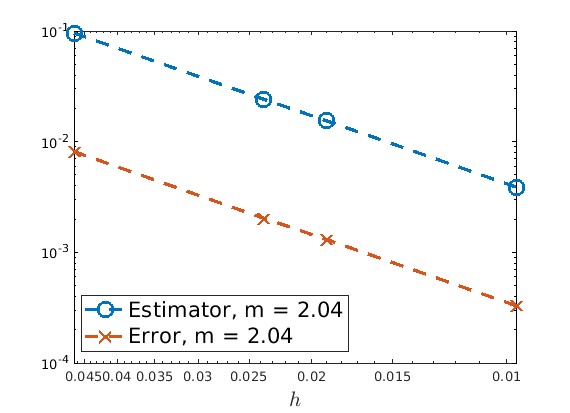}
    \caption{\meshtag{Distorted cartesian}, $\gamma_4$, $k=2$.}
    \label{subfig:DistortedCartesian_k2_Test3bis}
  \end{subfigure}\quad
     \begin{subfigure}[b]{.3\linewidth} \centering
    \includegraphics[width=1.\linewidth]{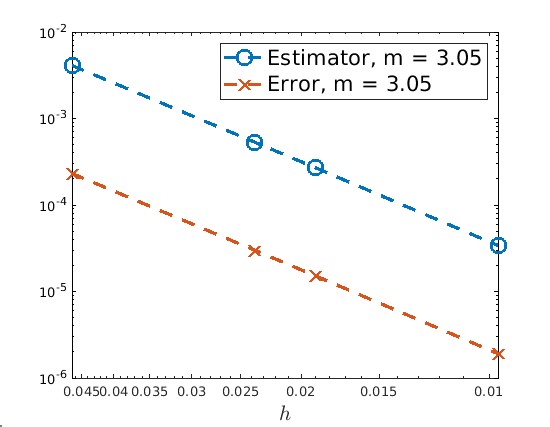}
    \caption{\meshtag{Distorted cartesian},  $\gamma_4$, $k=3$.}
    \label{subfig:DistortedCartesian_k3_Test3bis}
  \end{subfigure}
  \caption{Test 3: Convergence plots. In the legends, $m$ is the average convergence rate.}
  \label{Test3_convergence}
  \end{figure}

In this section, we consider problems where the diffusity coefficient does not meet the quasi-monotonicity condition discussed in Remark \ref{rmk:quasi-mon}. The tests considered are the one presented in \cite{Petzoldt2002}, where the diffusion tensor is $\K = \gamma_i I$ with $i =3,4$, with
\begin{align*}
\gamma_3(x,y) \coloneqq
\begin{cases}1
\quad\text{in }\Omega_{11} =  [0, 0.5)^2 ,\\
10^{-3} \quad\text{in } \Omega_{1 2} =[0.5, 1] \times[0, 0.5),\\
10^{-2} \quad\text{in } \Omega_{2 1} =[0, 0.5)\times[0.5, 1] ,\\
10 \quad\text{in } \Omega_{2 2} =[0.5, 1]^2 ,
\end{cases} &\gamma_4(x,y) \coloneqq
\begin{cases}1
\quad\text{in }\Omega_{11},\\
10^{-7} \quad\text{in } \Omega_{1 2},\\
10^{-2} \quad\text{in } \Omega_{2 1},\\
10^5 \quad\text{in } \Omega_{2 2} ,
\end{cases}
\end{align*}
We impose the loading terms in such a way that the exact solutions are
\begin{align*}
    u_i(x,y)= \xi_i(x) Y(y),
\end{align*}
where $\xi_i(x)$, $Y(y)$ are defined in \eqref{def:xi} and in \eqref{def:Yy}, and 
\begin{align*}
    c_i\coloneqq \begin{cases}
         - \frac{3 \gamma_i|_{\Omega_{1 1}}+ \gamma_i|_{\Omega_{1 2}}}{4\left(\gamma_i|_{\Omega_{1 1}}+ \gamma_i|_{\Omega_{1 2}}\right)}\quad \quad \text{in }\Omega_{11} \cup\Omega_{12},\\
         - \frac{3 \gamma_i|_{\Omega_{2 1}}+ \gamma_i|_{\Omega_{2 2}}}{4\left(\gamma_i|_{\Omega_{2 1}}+ \gamma_i|_{\Omega_{2 2}}\right)} \quad \quad \text{in }\Omega_{21} \cup\Omega_{22}.
    \end{cases}
\end{align*}
We solve the problem on a family of distorted cartesian meshes (see Figure \ref{subfig:DistortedCartesian_mesh}), conforming to the jumps of $\K$, achieving the optimal convergence rate as shown in Figure~\ref{Test3_convergence}.
Tables \ref{tab:distorted_Test3_gamma3}, and \ref{tab:distorted_Test3_gamma4} show the computed effectivity indices. These prove the robustness of the estimate proposed with respect to the meshsize. We notice also in this case that 
the jumps of $\K$ do not influence significantly  the effectivity indices, even though the quasi-monotonicity property is not fulfilled. Indeed, the area of the polygons where the quasi-monotonicity of $\K$ is not satisfied is decreasing with $h$.

\begin{table}[!ht]
    \centering
\begin{tabular}{c c |c c |c c}
\multicolumn{2}{c |}{$k$ = 1} &\multicolumn{2}{c |}{$k$ = 2} &\multicolumn{2}{c}{$k$ = 3} \\ 
\hline 
$h$ & $\epsilon$ &$h$ & $\epsilon$ &$h$ & $\epsilon$ \\ 
0.046649&6.2970&0.046649&11.7725&0.046649&17.7205 \\ 
0.023756&6.3248&0.023756&11.8629&0.023756&17.7857 \\ 
0.019008&6.3330&0.019008&11.8077&0.019008&17.7647 \\ 
0.0096754&6.3243&0.0096754&11.7982&0.0096754&17.8251 \\ 
\hline 
\end{tabular}
\caption{Test 3, case $\gamma_3$. \meshtag{Distorted cartesian}. Effectivity indices. }
    \label{tab:distorted_Test3_gamma3}
\end{table}

\begin{table}[!ht]
    \centering
\begin{tabular}{c c |c c |c c}
\multicolumn{2}{c |}{$k$ = 1} &\multicolumn{2}{c |}{$k$ = 2} &\multicolumn{2}{c}{$k$ = 3} \\ 
\hline 
$h$ & $\epsilon$ &$h$ & $\epsilon$ &$h$ & $\epsilon$ \\ 
0.046649&6.3282&0.046649&11.7832&0.046649&17.7266 \\ 
0.023756&6.3590&0.023756&11.8802&0.023756&17.7961 \\ 
0.019008&6.3682&0.019008&11.8213&0.019008&17.7691 \\ 
0.0096754&6.3567&0.0096754&11.8123&0.0096754&17.8354 \\ 
\hline 
\end{tabular}
\caption{Test 3, case $\gamma_4$. \meshtag{Distorted cartesian}. Effectivity indices. }
    \label{tab:distorted_Test3_gamma4}
\end{table}
\subsection{Test 4}
\begin{figure}[!ht]  
\centering 
\begin{subfigure}[b]{.3\linewidth} \centering
    \includegraphics[width=\linewidth]{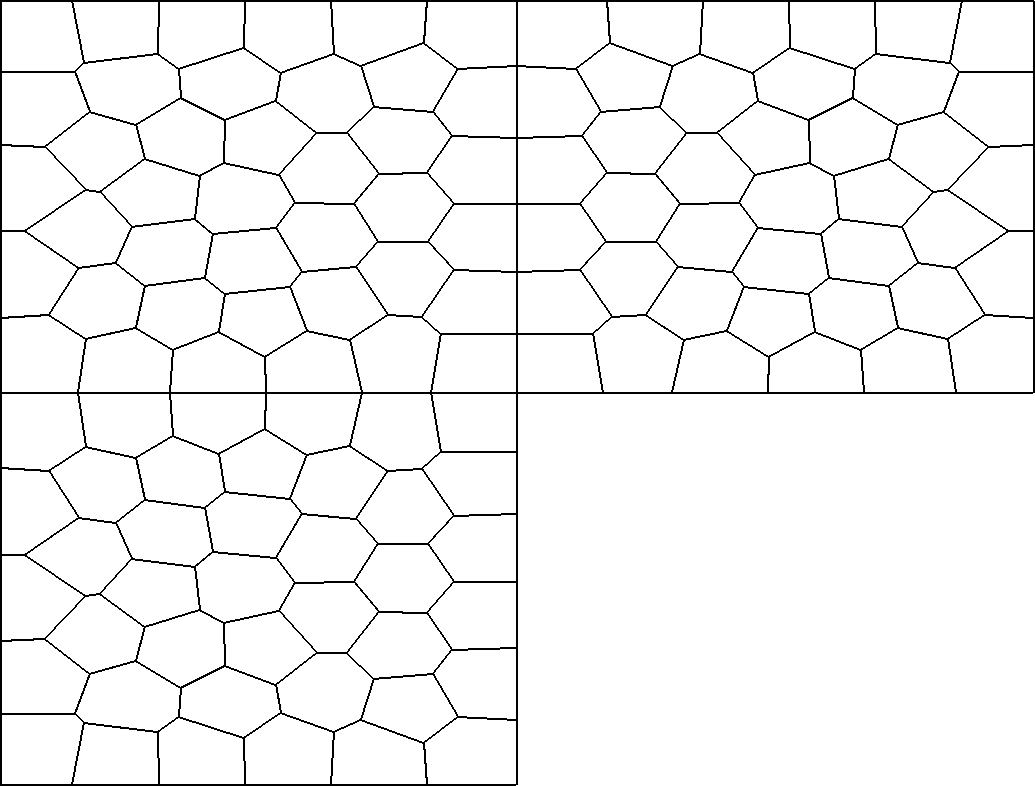}
    \caption{\meshtag{Polymesher}}
    \label{subfig:Lshape_Polymesher_mesh}
  \end{subfigure} \quad
\begin{subfigure}[b]{.3\linewidth} \centering
    \includegraphics[width=\linewidth]{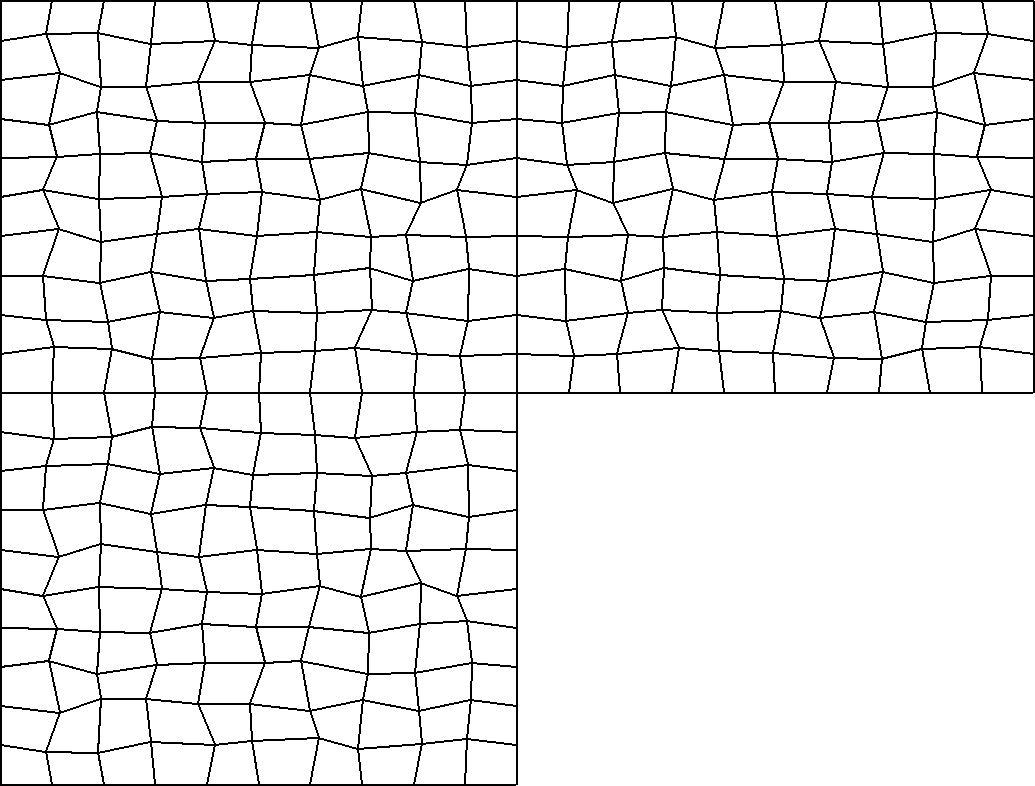}
    \caption{\meshtag{Distorted cartesian}}
    \label{subfig:Lshape_DistortedCartesian_mesh}
 \end{subfigure}  \quad
 \begin{subfigure}[b]{.3\linewidth} \centering
    \includegraphics[width=\linewidth]{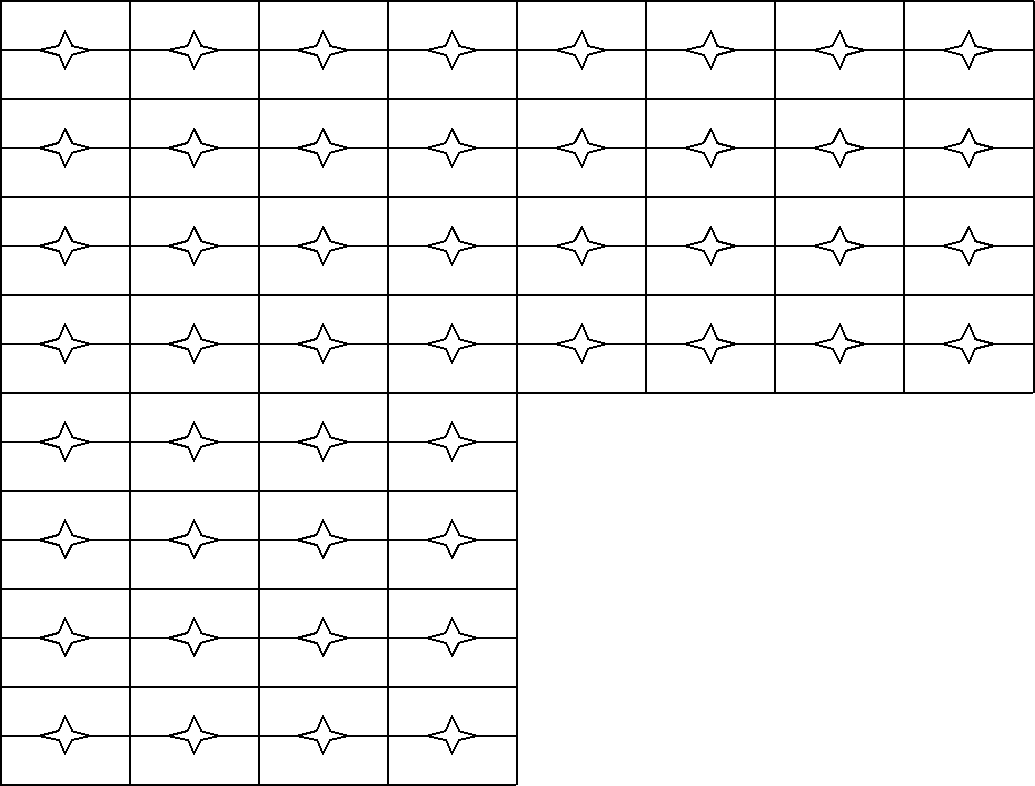}
    \caption{\meshtag{Star concave}}
    \label{Lshape_subfig:starconcave_mesh}
    \end{subfigure}
\caption{Meshes used for convergence Test 4.}
  \label{fig:Lshape_Meshes}
  \end{figure}

\begin{figure}
\begin{subfigure}[b]{.3\linewidth} \centering
    \includegraphics[width=1.\linewidth]{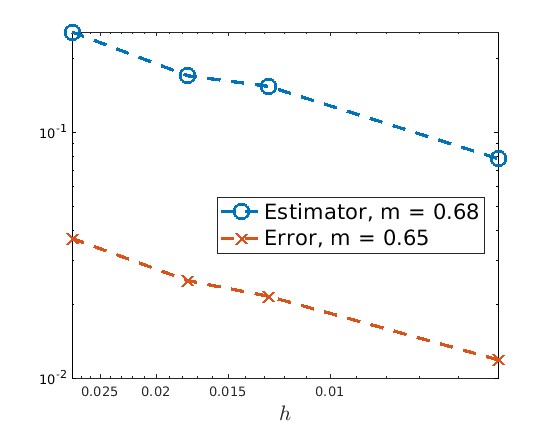}
    \caption{\meshtag{Polymesher}, $k=1$.}
    \label{subfig:Test4_Polymesher_k1}
  \end{subfigure}\quad
 \begin{subfigure}[b]{.3\linewidth} \centering
    \includegraphics[width=1.\linewidth]{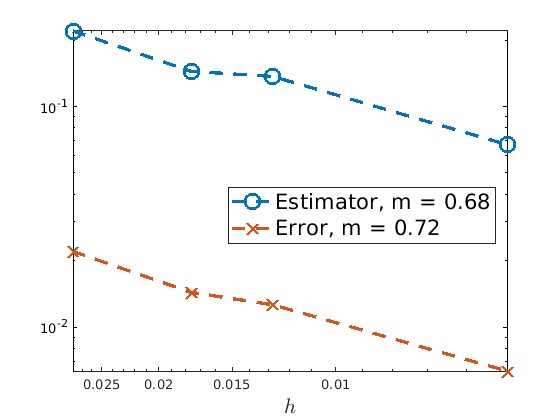}
    \caption{\meshtag{Polymesher}, $k=2$.}
    \label{subfig:Test4_Polymesher_k2}
  \end{subfigure}\quad
  \begin{subfigure}[b]{.3\linewidth} \centering
    \includegraphics[width=1.\linewidth]{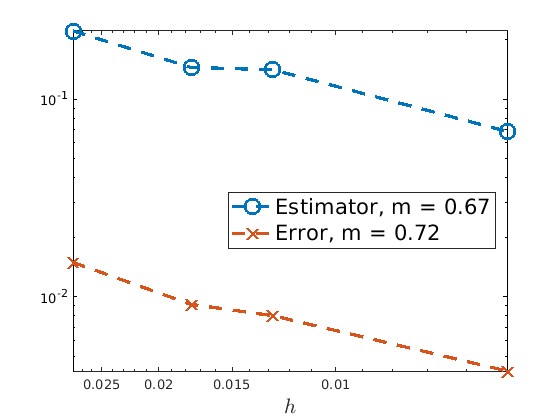}
    \caption{\meshtag{Polymesher}, $k=3$.}
    \label{subfig:Test4_Polymesher_k3}
  \end{subfigure}
  \\
  \begin{subfigure}[b]{.3\linewidth} \centering
    \includegraphics[width=1.\linewidth]{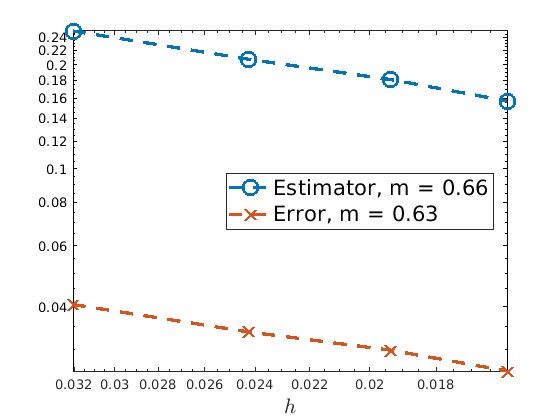}
    \caption{\meshtag{Distorted cartesian}, $k=1$.}
    \label{subfig:Test4_DistortedCartesian_k1}
  \end{subfigure}\quad
    \begin{subfigure}[b]{.3\linewidth} \centering
    \includegraphics[width=1.\linewidth]{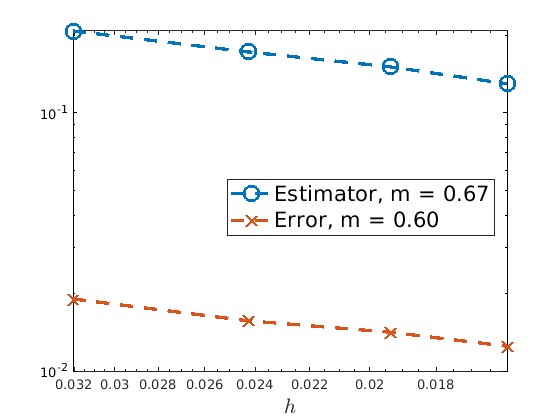}
    \caption{\meshtag{Distorted cartesian}, $k=2$.}
    \label{subfig:Test4_DistortedCartesian_k2}
  \end{subfigure}\quad
     \begin{subfigure}[b]{.3\linewidth} \centering
    \includegraphics[width=1.\linewidth]{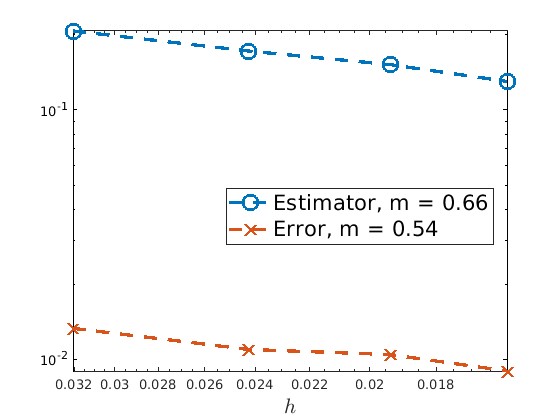}
    \caption{\meshtag{Distorted cartesian}, $k=3$.}
    \label{subfig:Test4_DistortedCartesian_k3}
  \end{subfigure}\\
    \begin{subfigure}[b]{.3\linewidth} \centering
    \includegraphics[width=1.\linewidth]{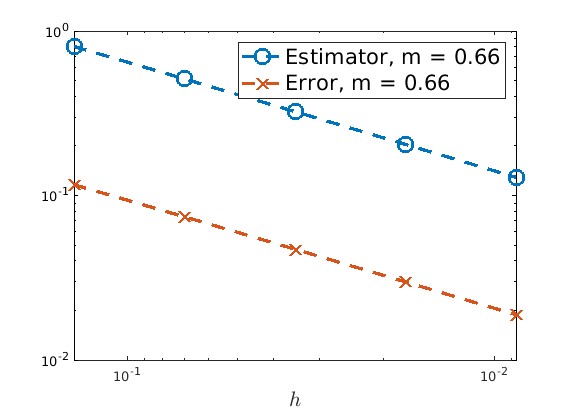}
    \caption{\meshtag{Star concave}, $k=1$.}
    \label{subfig:Test4_StarConcave_k1}
  \end{subfigure}\quad
    \begin{subfigure}[b]{.3\linewidth} \centering
    \includegraphics[width=1.\linewidth]{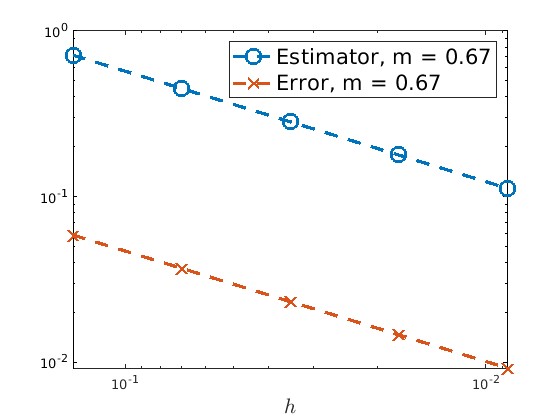}
    \caption{\meshtag{Star concave}, $k=2$.}
    \label{subfig:Test4_StarConcave_k2}
  \end{subfigure}\quad
     \begin{subfigure}[b]{.3\linewidth} \centering
    \includegraphics[width=1.\linewidth]{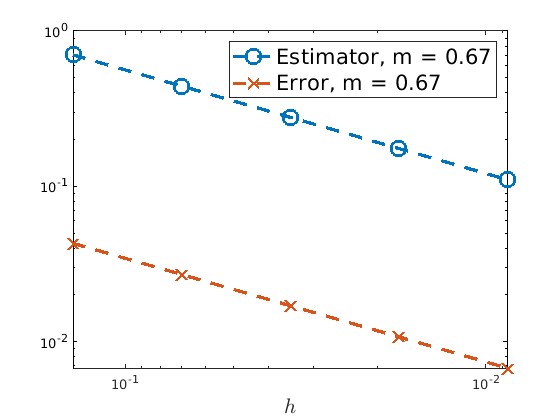}
    \caption{\meshtag{Star concave}, $k=3$.}
    \label{subfig:Test4_StarConcave_k3}
  \end{subfigure}
  \caption{Test 4: Convergence plots. In the legends, $m$ is the average convergence rate.}
  \label{Test4_convergence}
  \end{figure}

We consider a problem defined on the L-shaped domain $\Omega = (-1,1)^2\setminus \left([0,1]\times[-1,0] \right)$. We take the boundary conditions and the loading term such that the solution results
\begin{align*}
    u (\rho,\theta) = \rho^{2/3},
\end{align*}
where $\rho$ and $\theta$ are the polar coordinates. This solution $u\in\sobh{1}{\Omega}$ has a corner singularity in the origin, and it is possible to prove that $u\in \sobh{s}{\Omega}$, $s<\frac{5}{3}$. Under uniform mesh refinements,  the asymptotic rate is suboptimal and it results $\left(\# \Mh \right)^{-1/3} \simeq h^{2/3}$, independently on the polynomial degree $k$. We performed the test on \meshtag{Polymesher}, \meshtag{Distorted cartesian}, and \meshtag{Star concave} L-shaped meshes (see Figure \ref{fig:Lshape_Meshes}) and for the cases $k=1,2,3$, achieving the expected order of convergence as shown in Figure \ref{Test4_convergence}. The test highlights that the effectivity indices are not affected by the bounded Sobolev regularity of the problem, as shown in Tables~\ref{tab:Polymesher_Test4}, \ref{tab:distorted_Test4}, \ref{tab:StarConcave_Test4}.

\begin{table}[!ht]
   \centering
\begin{tabular}{c c |c c |c c}
\multicolumn{2}{c |}{$k$ = 1} &\multicolumn{2}{c |}{$k$ = 2} &\multicolumn{2}{c}{$k$ = 3} \\ 
\hline 
$h$ & $\epsilon$ &$h$ & $\epsilon$ &$h$ & $\epsilon$ \\ 
0.027988&6.9189&0.027988&10.0000&0.027988&14.8305 \\ 
0.017629&6.8200&0.017629&10.0603&0.017629&15.8111 \\ 
0.012807&7.1462&0.012807&10.8552&0.012807&17.5659 \\ 
0.005119&6.5819&0.005119&10.6305&0.005119&16.3231 \\ 
\hline 
\end{tabular}
\caption{Test 4, \meshtag{Polymesher}. Effectivity indices.}
    \label{tab:Polymesher_Test4}
\end{table}

\begin{table}[!ht]
    \centering
\begin{tabular}{c c |c c |c c}
\multicolumn{2}{c |}{$k$ = 1} &\multicolumn{2}{c |}{$k$ = 2} &\multicolumn{2}{c}{$k$ = 3} \\ 
\hline 
$h$ & $\epsilon$ &$h$ & $\epsilon$ &$h$ & $\epsilon$ \\ 
0.032068&6.1629&0.032068&10.9304&0.032068&15.5400 \\ 
0.024269&6.1234&0.024269&11.0116&0.024269&15.7007 \\ 
0.019351&6.0442&0.019351&10.6678&0.019351&14.5184 \\ 
0.016079&6.0182&0.016079&10.4016&0.016079&14.4664 \\ 
\hline 
\end{tabular}
\caption{Test 4, \meshtag{Distorted cartesian}. Effectivity indices.}
    \label{tab:distorted_Test4}
\end{table}

\begin{table}[!ht]
    \centering
\begin{tabular}{c c |c c |c c}
\multicolumn{2}{c |}{$k$ = 1} &\multicolumn{2}{c |}{$k$ = 2} &\multicolumn{2}{c}{$k$ = 3} \\ 
\hline 
$h$ & $\epsilon$ &$h$ & $\epsilon$ &$h$ & $\epsilon$ \\ 
0.13975&6.9244&0.13975&12.1802&0.13975&16.3426 \\ 
0.069877&6.8886&0.069877&12.1814&0.069877&16.3432 \\ 
0.034939&6.8648&0.034939&12.1820&0.034939&16.3434 \\ 
0.017469&6.8493&0.017469&12.1822&0.017469&16.3435 \\ 
0.0087346&6.8395&0.0087346&12.1822&0.0087346&16.3436 \\ 
\hline 
\end{tabular}
\caption{Test 4, \meshtag{Star concave}. Effectivity indices.}
    \label{tab:StarConcave_Test4}
\end{table}

\section{Conclusions}
\label{sec:conclusions}

We derived \emph{a posteriori} error estimates for the Stabilization-Free
Virtual Element Method defined in \cite{Berrone2024a,BBFM2025}. The abscence of
a stabilizing bilinear form in the discretization method allows to obtain
equivalence between a suitably defined error measure and classical residual
error estimators. Numerical tests on several mesh types validate the proposed
error estimator.

\section*{Acknowledgments}
The four authors are members of the Gruppo Nazionale Calcolo Scientifico (GNCS) at Istituto
Nazionale di Alta Matematica (INdAM). 
The authors kindly acknowledge financial support by INdAM-GNCS, by the Italian Ministry of University and Research (MUR) and by the European Union.
\section*{Fundings}
The four authors are members of the Gruppo Nazionale Calcolo Scientifico (GNCS) at Istituto
Nazionale di Alta Matematica (INdAM). The authors kindly acknowledge financial support by INdAM-GNCS
Project 2025 (CUP: E53C24001950001), by the Italian Ministry of University and Research (MUR) through 
the project MUR-M4C2-1.1-PRIN 2022\\ (CUP: E53D23005820006), and by the European Union through Next
Generation EU, M4C2, PRIN 2022 PNRR project (CUP: E53D23017950001) and PNRR M4C2 project of CN00000013 National Centre for HPC, Big Data and Quantum
Computing (HPC) (CUP: E13C22000990001).
\bibliographystyle{plain}
\bibliography{bibliografia}
\end{document}